\documentclass[UTF-8,reqno]{amsart}
\usepackage{enumerate}
\usepackage{amssymb,url,color, booktabs}
\usepackage{mathrsfs}

\usepackage{color}
\definecolor{MyDarkBlue}{cmyk}{0.8,0.3,0.8,0.4}
\definecolor{yellow}{rgb}{0.99,0.99,0.70}
\definecolor{white}{rgb}{1.0,1.0,1.0}
\definecolor{black}{rgb}{0.00,0.00,0.00}

\numberwithin{equation}{section}

\newcommand{\be}{\begin{eqnarray}}
\newcommand{\ee}{\end{eqnarray}}
\newcommand{\ce}{\begin{eqnarray*}}
\newcommand{\de}{\end{eqnarray*}}
\newtheorem{theorem}{Theorem}[section]
\newtheorem{lemma}[theorem]{Lemma}
\newtheorem{remark}[theorem]{Remark}
\newtheorem{definition}[theorem]{Definition}
\newtheorem{proposition}[theorem]{Proposition}
\newtheorem{Examples}[theorem]{Example}
\newtheorem{corollary}[theorem]{Corollary}

\def\nor{|\mspace{-3mu}|\mspace{-3mu}|}

\def\eps{\varepsilon}

\def\u{\mathbf{u}}

\def\p{\partial}

\def\[{{\Big[}}
\def\]{{\Big]}}
\def\<{{\langle}}
\def\>{{\rangle}}
\def\({{\Big(}}
\def\){{\Big)}}

\def\bx{{\mathbf{x}}}

\def\dif{{\mathord{{\rm d}}}}

\def\min{{\mathord{{\rm min}}}}

\def\no{\nonumber}
\def\={&\!\!=\!\!&}

\def\cB{{\mathcal B}}

\def\cI{{\mathcal I}}
\def\cJ{{\mathcal J}}
\def\cK{{\mathcal K}}

\def\cM{{\mathcal M}}

\def\cP{{\mathcal P}}

\def\mC{{\mathbb C}}
\def\mD{{\mathbb D}}
\def\mE{{\mathbb E}}

\def\mH{{\mathbb H}}
\def\mI{{\mathbb I}}

\def\mL{{\mathbb L}}

\def\mN{{\mathbb N}}

\def\mP{{\mathbb P}}
\def\mQ{{\mathbb Q}}
\def\mR{{\mathbb R}}
\def\mS{{\mathbb S}}

\def\mZ{{\mathbb Z}}

\def\bP{{\mathbf P}}

\def\bE{{\mathbf E}}
\def\1{{\mathbf{1}}}

\def\sC{{\mathscr C}}
\def\sD{{\mathscr D}}
\def\sE{{\mathscr E}}
\def\sF{{\mathscr F}}
\def\sG{{\mathscr G}}

\def\sI{{\mathscr I}}
\def\sJ{{\mathscr J}}

\def\sL{{\mathscr L}}

\def\sV{{\mathscr V}}

\def\geq{\geqslant}
\def\leq{\leqslant}

\def\div{\mathord{{\rm div}}}

\def\eps{\varepsilon}

\def\u{\mathbf{u}}

\def\p{\partial}

\def\[{{\Big[}}
\def\]{{\Big]}}
\def\<{{\langle}}
\def\>{{\rangle}}
\def\({{\Big(}}
\def\){{\Big)}}

\def\bx{{\mathbf{x}}}

\def\dif{{\mathord{{\rm d}}}}

\def\min{{\mathord{{\rm min}}}}

\def\no{\nonumber}
\def\={&\!\!=\!\!&}
\def\bt{\begin{theorem}}
\def\et{\end{theorem}}
\def\bl{\begin{lemma}}
\def\el{\end{lemma}}
\def\br{\begin{remark}}
\def\er{\end{remark}}
\def\bx{\begin{Examples}}
\def\ex{\end{Examples}}
\def\bd{\begin{definition}}
\def\ed{\end{definition}}
\def\bp{\begin{proposition}}
\def\ep{\end{proposition}}
\def\bc{\begin{corollary}}
\def\ec{\end{corollary}}

\def\geq{\geqslant}
\def\leq{\leqslant}

\def\div{\mathord{{\rm div}}}

\def\bP{{\mathbf P}}

\def\<{\langle} \def\>{\rangle}

\allowdisplaybreaks

\begin{document}


\title[Maximum principle for non-uniformly parabolic equations]
{Maximum principle  for non-uniformly parabolic equations and applications}

\author{Xicheng Zhang}

\address{Xicheng Zhang:
School of Mathematics and Statistics, Wuhan University,
Wuhan, Hubei 430072, P.R.China\\
Email: XichengZhang@gmail.com
 }

\thanks{
This work is partially supported by NNSFC grant of China (No. 11731009), and the German Research Foundation (DFG) through the Collaborative Research Centre(CRC) 1283 ``Taming uncertainty and profiting from randomness and low regularity in analysis, stochastics and their applications".
}

\begin{abstract}
In this paper we study the global boundedness for the solutions to a class of possibly degenerate parabolic equations by De-Giorgi's iteration.
As applications, we show the existence of weak solutions for possibly degenerate stochastic differential equations with 
singular diffusion and drift coefficients. Moreover, by the Markov selection theorem of Krylov [8], 
we also establish the existence of the associated strong Markov family.


\end{abstract}

\maketitle \rm

\section{Introduction}

Consider the following elliptic equation of divergence form in $\mR^d$ ($d\geq 2$):
\begin{align}\label{ELL}
\div (a\cdot\nabla u)=0,
\end{align}
where $a:\mR^d\to \mR^{d\times d}$ is a Borel measurable function and $\nabla:=(\p_{x_1},\cdots,\p_{x_d})$. We introduce the following two functions:
\begin{align}\label{DZ1}
\lambda_0(x):=\inf_{|\xi|=1}\xi\cdot a(x)\xi,\quad \mu_0(x):=\sup_{|\xi|=1}\frac{|a(x)\xi|^2}{\xi\cdot a(x)\xi}.
\end{align}
Suppose that $\lambda_0$ and $\mu_0$ are nonnegative measurable functions.
If $\lambda^{-1}_0$ and $\mu_0$ are essentially bounded, that is, $a$ is uniformly elliptic, then the celebrated works of De-Giorgi \cite{De-Gi} and Nash \cite{Na} said that
any weak solutions of elliptic equation \eqref{ELL} are bounded and H\"older continuous. Moreover, Moser \cite{Mo} showed that any weak solutions of \eqref{ELL}
satisfy the Harnack inequality.

\medskip

In \cite{Tr}, Trudinger considered the non-uniformly elliptic equation \eqref{ELL} under the following integrability assumptions:
$$
\lambda^{-1}_0\in L^{p_0},\ \mu_0\in L^{p_1}\mbox{ with $p_0,p_1\in(1,\infty]$ satisfying $\tfrac{1}{p_0}+\tfrac{1}{p_1}<\tfrac{2}{d}$,}
$$
and showed that any generalized solutions of \eqref{ELL} are locally bounded and weak Harnack inequality holds.
Recently, Bella and Sch\"affner \cite{Be-Sc} showed the same results under the following {\it sharp} condition on $p_0, p_1$,
\begin{align}\label{PP1}
\tfrac{1}{p_0}+\tfrac{1}{p_1}<\tfrac{2}{d-1},\ \ p_0,p_1\in[1,\infty],
\end{align}
where the key point is a new Sobolev embedding inequality of  variational type. In this paper we are interested in a parabolic version of \cite{Be-Sc}, 
and aim to establish the global boundedness for the solutions of non-uniformly parabolic equations. 
More precisely, we shall consider the following parabolic equation of divergence form in $\mR^{d+1}$:
\begin{align}\label{PDE0}
\p_t u=\div (a\cdot\nabla u)+b\cdot \nabla u+f,
\end{align}
where 
$$
a:\mR^{d+1}\to\mR^{d\times d},\ b:\mR^{d+1}\to\mR^d,\ f:\mR^{d+1}\to\mR
$$
are Borel measurable functions. As in \eqref{DZ1}, we introduce 
\begin{align}\label{DW3}
\lambda(x):=\inf_{t\geq 0, |\xi|=1}\xi\cdot a(t,x)\xi,\quad \mu(x):=\sup_{t\geq 0, |\xi|=1}\frac{|a(t,x)\xi|^2}{\xi\cdot a(t,x)\xi},
\end{align}
and suppose that $\lambda$ and $\mu$ are nonnegative Borel measurable functions.

\smallskip

First of all we introduce the following notion of weak solutions to PDE \eqref{PDE0}.
\bd\label{Def1}
A continuous function $u:\mR^{d+1}\to\mR$ is called a Lipschitz weak (super/sub)-solution of PDE \eqref{PDE0} 
if $\nabla u$ is locally bounded and for any nonnegative Lipschitz function $\varphi$ on $\mR^{d+1}$ with compact support,
\begin{align}\label{Def0}
-\<\!\<u,\p_t\varphi\>\!\>=(\geq/\leq)-\<\!\<a\cdot\nabla u,\nabla\varphi\>\!\!+\<\!\<b\cdot\nabla u,\varphi\>\!\>+\<\!\<f,\varphi\>\!\>,
\end{align}
where $\<\!\<f,g\>\!\>:=\int_\mR\int_{\mR^{d}}f(t,x)g(t,x)\dif x\dif t$.
\ed
Throughout this paper, we fix $p_0\in(\frac{d}{2},\infty]$ and $p_1\in[1,\infty]$ with
\begin{align}\label{PP0}
\tfrac1{p_0}+\tfrac{1}{p_1}<\tfrac{2}{d-1}.
\end{align}
With the notations in \eqref{FQ22} and \eqref{FQ2} below,
we assume that
\begin{enumerate}[{\bf (H$^a$)}]
\item  $\nor\lambda^{-1}\nor_{p_0}+\nor\mu\nor_{p_1}<\infty$, where $\lambda,\mu$ are defined by \eqref{DW3}.
\end{enumerate}
\begin{enumerate}[{\bf (H$^b$)}]
\item $b=b_1+b_2$, where 
if $p_0\in(\tfrac d2,d]$, $b_1\equiv 0$,  and  if $p_0>d$, $b_1\in\widetilde\mL^{q_2,p_2}_{t,x}$ for some $(p_2,q_2)\in[1,\infty]^2$ with
\begin{align}\label{Re1}
\tfrac{1}{2p_0}+\tfrac{1}{p_2}<(\tfrac12-\tfrac{1}{q_2})(\tfrac{2}{d}-\tfrac{1}{p_0}),
\end{align}
and $(\div b_2)^-=0$ and $b_2\in\widetilde\mL^{p_3,q_3}_{x,t}$ for some $1\leq p_3\leq q_3\leq\infty$ with
\begin{align}\label{Re01}
\tfrac{(d-1)\vartheta_1}{p_3}+\tfrac{2+\vartheta_1+d(\vartheta_2-\vartheta_1)}{q_3}<2,
\end{align}
where $\vartheta_1:=\big(1-\tfrac{d-1}{2p_0}\big)^{-1}$ and $\vartheta_2:=\big(1-\tfrac{d}{2p_0}\big)^{-1}$.
\end{enumerate}
\br\rm
Note that condition \eqref{Re1} is satisfied for $p_2=q_2=\infty$ if and only if $p_0>d$. This is why we need to put $b_1\equiv 0$ for $p_0\leq d$.
If $p_0=\infty$, i.e., $a$ has a lower bound, condition \eqref{Re1} reduces to the usual one $\frac{d}{p_2}+\frac{2}{q_2}<1$,
and condition \eqref{Re01} becomes $\frac{d-1}{p_3}+\frac{3}{q_3}<2$ (see Corollary \ref{Cor1} below).
\er

For simplicity of notations, we introduce the parameter set
\begin{align}\label{TH}
\Theta:=\big(d,p_i, q_i,\nor\lambda^{-1}\nor_{{p_0}},\nor\mu\nor_{p_1}, \nor b_1\nor_{\widetilde\mL^{q_2,p_2}_{t,x}}, \nor b_2\nor_{\widetilde\mL^{p_3,q_3}_{x,t}}\big),
\end{align}
and the index set
$$
\mI^d_{p_0}:=\Big\{(p,q)\in[1,\infty]^2: \tfrac{1}{p}<(1-\tfrac{1}{q})(\tfrac{2}{d}-\tfrac{1}{p_0})\Big\}.
$$
The main aim of this paper is to prove the following a priori estimate.
\bt\label{TH22}
Under {\bf (H$^a$)} and {\bf (H$^b$)}, for any $f\in\widetilde\mL^{q_4,p_4}_{t,x}$ with $(p_4,q_4)\in\mI^d_{p_0}$ and for  any $T>0$,
there exists a constant $C=C(T,\Theta, p_4,q_4)>0$ such that for any Lipschitz weak solution $u$ of PDE \eqref{PDE0} in $\mR^{d+1}$ 
with $u(t)|_{t\leq 0}\equiv 0$,
\begin{align}\label{Es161}
\|u\|_{L^\infty([0,T]\times\mR^d)}+\nor u \1_{[0,T]}\nor_{\widetilde\sV}
\leq C\nor f\1_{[0,T]}\nor_{\widetilde\mL^{q_4,p_4}_{t,x}},
\end{align}
where $\widetilde\sV$ is defined by $\eqref{VV}$ below.
\et
\br\rm
After this paper was posted on the arXiv, I learned that Bella and Sch\"affner \cite{Be-Sc1} recently proved a non-uniformly parabolic version of their results 
for finite-difference operators of divergence form in lattice $\mZ^d$. 
I would like to mention that to overcom the difficulty caused by the time variable, we use completely different methods
(see \cite[Lemma 2]{Be-Sc1} and Lemma 3.2 below). Moreover, the novelty of this paper is that we are considering the supercritical drift $b$ and nonhomogeneous $f$, 
which are crucial for applications in SDEs.
\er

Consider the following heat equation with divergence free drift $b$:
\begin{align}\label{DQ1}
\p_t u=\Delta u+b\cdot \nabla u+f,\ u(t)|_{t\leq 0}=0.
\end{align}
By Theorem \ref{TH22} with $p_0=\infty$,
we obtain the following a priori global boundedness estimate, to the author's knowledge, which is new. 
\bc\label{Cor1}
Let $b\in\widetilde\mL^{p,q}_{x,t}$ with $\div b=0$, where $1\leq p\leq q\leq\infty$ satisfy
\begin{align}\label{PP9}
\tfrac{d-1}p+\tfrac3q<2.
\end{align}
For any $T>0$ and $f\in\widetilde\mL^{q',p'}_{t,x}$, where $p',q'\in[1,\infty]$ satisfy $\frac{d}{p'}+\frac{2}{q'}<2$,
there exists a constant $C>0$ only depending on $T,d,p,p',q'$ and $\|b\|_{\widetilde\mL^{p,q}_{x,t}}$ such that for any Lipschitz weak solution $u$ of \eqref{DQ1},
\begin{align}\label{Es162}
\|u\|_{L^\infty([0,T]\times\mR^d)}\leq C\nor f\1_{[0,T]}\nor_{\widetilde\mL^{q',p'}_{t,x}}.
\end{align}
\ec
\br\rm
For divergence free drift $b$, when $b\in\widetilde\mL^{q,p}_{t,x}$ for some $\tfrac{d}{p}+\tfrac{2}{q}<2$, it is well known that \eqref{Es162}
holds (cf. \cite{Na-Ur}, \cite{Zh-Zh}). However, when $b$ does not depend on the time variable $t$, the condition \eqref{PP9}
in Corollary \ref{Cor1} is clearly better than $p>\frac d2$. Notice that for $p=q$, they are same. 
In a forthcoming paper, we shall establish the Harnack inequality for PDE \eqref{DQ1} under \eqref{PP9} (see \cite{Ig-Ku-Ry}).
Here, an interesting open question is that for divergence free and time independent drift $b(x)$, whether condition \eqref{PP9} is sharp for \eqref{Es162}.
\er

In \cite{Be-Sc}, the local boundedness of generalized solutions of elliptic equations 
is used to establish the $L^\infty$-sublinearity of the corrector in stochastic homogenization
 in non-uniformly case, which is a key step of proving quenched invariance principle for random walks \cite{An-De-Sl}. 
 As in \cite{Be-Sc} and \cite{An-Ch-De-Sl}, Theorem \ref{TH22} could be used to showing a quenched invariance principle for random walks 
 in time-dependent ergodic environment (see \cite{Be-Sc1} for recent development). 
As one application of the global boundedness estimate \eqref{Es161}, 
we shall establish the existence of weak solutions to possibly degenerate SDEs with singular diffusion and drift coefficients in this paper. 
Consider the following SDE:
\begin{align}\label{SDE0}
\dif X_t=\sqrt{2}\sigma(t,X_t)\dif W_t+b(t,X_t)\dif t,\ \ X_0=x,
\end{align}
where $W$ is a $d$-dimensional standard Brownian motion on some stochastic basis $(\Omega,\sF,\mP; (\sF_t)_{t\geq 0})$
and $\sigma:\mR_+\times\mR^d\to\mR^d\otimes\mR^d$
and $b:\mR_+\times\mR^d\to\mR^d$ are Borel measurable functions. Note that the generator of SDE \eqref{SDE0} is given by
$$
\sL^{\sigma,b}_t f(x)=(\sigma^{ik}\sigma^{jk})(t,x)\p_i\p_j f(x)+b^j(t,x)\p_j f(x).
$$
Here and after we shall use the usual Einstein convention for summation: an index appearing twice in a product will be summed automatically.

\medskip

It is well known that if $\sigma$ and $b$ are Lipschitz continuous in $x$ uniformly in $t$, then SDE \eqref{SDE0} admits a unique strong solution.
When  $\sigma$ is bounded measurable and uniformly elliptic and $b\in L^{d+1}(\mR_+\times\mR^d)$, recently, Krylov \cite{Kr20} 
showed the existence of weak solutions to  SDE \eqref{SDE0} (see \cite{Kr80} for bounded measurable drift $b$).
When $\sigma$ is the identity matrix and $b$ is divergence free and belongs to $\widetilde\mL^{q,p}_{t,x}$ for some $p,q\in[1,\infty]$
with $\frac{d}{p}+\frac{2}{q}<2$, utilizing the like-estimate \eqref{Es162}, in a joint work \cite{Zh-Zh} with G. Zhao, we showed the existence of weak solutions to SDE \eqref{SDE0}.
In particular, the stochastic Lagrangian trajectories associated with Leray's solutions of 3D-Navier-Stokes equations are constructed.
However, when diffusion coefficient $\sigma$ is possibly degenerate or singular, and $b$ is irregular (saying only bounded measurable), 
to the author's knowledge, it seems that there are few results about the existence of solutions to SDE \eqref{SDE0} except for \cite{Wa-Zh}. 
To show the existence of weak solutions, the key step is to prove the following estimate of  Krylov's type: for any $(p,q)\in\mI^d_{p_0}$,
\begin{align}\label{Kry10}
\mE\left(\int^t_0 f(s, X_s)\dif s\right)\leq C\nor f\nor_{\widetilde\mL^{q,p}_{t,x}}.
\end{align}
Note that if we let $a=\sigma\sigma^*$, then $\sL^{\sigma,b}_t$ can be written as the divergence form:
$$
\sL^{\sigma,b}_t f(x)=\p_i(a^{ij}(t,\cdot)\p_j f)(x)+(b^j-\p_ia^{ij})(t,x)\p_j f(t,x).
$$
Under suitable conditions, \eqref{Kry10} will be a consequence of It\^o's formula and \eqref{Es161} (see Theorem \ref{Th43} below). 

\medskip

Although we can show the existence of weak solutions for SDE \eqref{SDE0} with singular coefficients, 
in many cases, the uniqueness is not easily obtained and even does not hold for SDEs with measurable coefficients. In 1973,
N.V. Krylov \cite{Kr73} proved a Markov selection theorem from the family of solutions of SDE \eqref{SDE0} 
when $b$ and $\sigma$ are bounded continuous.
His method was presented in a different way in \cite[Chapter 12]{St-Va}. For applications in SPDEs, we refer to \cite{Fl-Ro} and \cite{Go-Ro-Zh}.
Here we shall follow Stroock and Varadhan's method \cite{St-Va} to select a strong Markovian solution 
for SDEs \eqref{SDE0} with singular coefficients when the uniqueness is not applicable. 

\medskip

We would like to mention the following examples to illustrate our main results obtained in Sections 4 and 5.
\bx\rm
Let $d=3$ and $\u(t,x)$ be any Leray solutions of 3D-Navier-Stokes equations. Consider the following SDEs:
$$
\dif X_{t,s}=\sqrt{2}\dif W_t+\u(t,X_{t,s})\dif t,\ t\geq s\geq 0, \ X_{s,s}=x\in\mR^3.
$$
In \cite{Zh-Zh}, the existence of weak solutions is obtained to the above SDE. 
By \cite[Theorem 1.1]{Zh-Zh} and Theorem \ref{Th55} below, one can select a family of probability measures
$(\mP_{s,x})_{(s,x)\in\mR_+\times\mR^3}$ on the continuous function space $\mC$ so that 
for each $(s,x)\in\mR_+\times\mR^3$, $\mP_{s,x}$ solves the martingale problem associated to the above SDE, and
$(\mP_{s,x})_{(s,x)\in\mR_+\times\mR^3}$ forms a time-inhomoegenous 
strong Markovian family. 
\ex

\bx\rm
Let $d\geq 3$ and $\alpha\in(0,(\frac d2-1)\wedge(\frac12+\frac{1}{d-1}))$, $\beta\in(0,2\alpha)$. For any $\lambda\geq0$ and $x\in\mR^d$, 
the following SDE admits a unique strong solution (see Proposition \ref{Prop62} below):
$$
\dif X_t=|X_t|^{-\alpha}\dif W_t+\lambda X_t|X_t|^{-\beta-1}\dif t,\ \  X_0=x.
$$
Note that the starting point can be zero.
\ex

\smallskip

This paper is organized as follows: In Section 2, we prove a time-dependent variational embedding theorem, which in particular extends the result obtained in \cite{Be-Sc}.
In Section 3, we prove our main Theorem \ref{TH22} by De-Giorgi's iteration (cf. \cite{De-Gi}). 
In Section 4, we apply our main result to SDEs with rough coefficients. In Section 5, we use Krylov's Markov selection theorem
to select a strong Markov family from the weak solution family. In Section 6, we present two examples to illustrate our result.
In the appendix, we recall some results about the
regular conditional probability distribution (abbreviated as r.c.p.d.) as well as the abstract time-inhomoegenous 
strong Markov selection theorem.

\medskip

Throughout this paper, we use the following conventions: The letter $C=C(\cdots)$ denotes a constant, whose value may change in different places, and which is increasing with 
respect to its argument.
We also use $A\lesssim B$  or $A\lesssim_CB$ to denote $A\leq C B$ for some unimportant constant $C>0$. 

\section{Preliminaries}

Let $\sD:=C^\infty_c(\mR^{d+1})$ 
be the space of all smooth functions in $\mR^{d+1}$ with compact supports and $\sD'$ the dual space of 
$\sD$, which is also called the distribution space. The duality between $\sD'$ and $\sD$ is denoted by $\<\!\<\cdot,\cdot\>\!\>$.
In particular, if $f\in\sD'$ is locally integrable and $g\in\sD$, then
\begin{align}\label{EP1}
\<\!\<f,g\>\!\>=\int_\mR\<f(t),g(t)\>\dif t\ \ \mbox{with\ \ } \<f(t),g(t)\>:=\int_{\mR^d}f(t,x)g(t,x)\dif x.
\end{align}
For $p,q\in[1,\infty]$, let $\mL^{q,p}_{t,x}:=L^q(\mR; L^p(\mR^d))$
and $\mL^{p,q}_{x,t}:=L^p(\mR^d; L^q(\mR))$ be the space of spatial-time functions with norms, respectively,
$$
\|f\|_{\mL^{q,p}_{t,x}}:=\left(\int_{\mR}\|f(t,\cdot)\|_p^q\dif t\right)^{1/q},\ 
\|f\|_{\mL^{p,q}_{x,t}}:=\left(\int_{\mR^d}\|f(\cdot,x)\|_q^p\dif x\right)^{1/p},
$$
where $\|\cdot\|_p$ stands for the usual $L^p$-norm.
By Minkowskii's inequality,
\begin{align}\label{HQ5}
\|f\|_{\mL^{q,p}_{t,x}}\leq \|f\|_{\mL^{p,q}_{x,t}}\mbox{ if } q\geq p;\quad \|f\|_{\mL^{p,q}_{x,t}}\leq \|f\|_{\mL^{q,p}_{t,x}}\mbox{ if } q\leq p.
\end{align}
For $r>0$ and $(s,z)\in\mR^{d+1}$, we define
$$
Q_r:=[-r^2,r^2]\times B_r\subset\mR^{d+1},\ Q^{s,z}_r:=Q_r+(s,z),\ \ B^z_r:=B_r+z,
$$
and for $p\in[1,\infty]$, introduce the following localized $L^p$-space:
\begin{align}\label{FQ22}
{\widetilde L}^p:=\Big\{f\in L^1_{loc}(\mR^d): \nor f\nor_p:=\sup_z\|\1_{B^z_1}f\|_p<\infty\Big\},
\end{align}
and for $p,q\in[1,\infty]$,
\begin{align}\label{FQ2}
\widetilde\mL^{q,p}_{t,x}:=\Big\{f\in L^1_{loc}(\mR^{d+1}): \nor f\nor_{\widetilde\mL^{q,p}_{t,x}}:=\sup_{s,z}\|\1_{Q^{s,z}_1}f\|_{\mL^{q,p}_{t,x}}<\infty\Big\},
\end{align}
and similarly for $\widetilde\mL^{p,q}_{x,t}$. 
Clearly, for $p\leq p'$ and $q\leq q'$, 
$$
\widetilde\mL^{q',p'}_{t,x}\subset \widetilde\mL^{q,p}_{t,x},\ \ \widetilde\mL^{p',q'}_{x,t}\subset \widetilde\mL^{p,q}_{x,t}.
$$
By a finite covering technique, it is easy to see that for any $T,r>0$ (see \cite{Zh-Zh}),
\begin{align}\label{DK1}
\nor \1_{[0,T]}f\nor_{\widetilde\mL^{q,p}_{t,x}}\asymp\sup_{z}\|\1_{[0,T]\times B^z_r}f\|_{\mL^{q,p}_{t,x}}.
\end{align}

First of all, we have the following Gagliado-Nirenberge's interpolation estimate.
\bl\label{Le11}
Fix ${\varkappa}\in[2d/(d+2),2]$ and $\theta\in[0,1]$ with exception $\theta=1$ and ${\varkappa}=d$. For any $r\geq 2$ and $s\geq 1$ with
$$
\tfrac{1}{2}-\tfrac{1}{r}=\tfrac\theta 2\Big(\tfrac{2}{d}+1-\tfrac{2}{{\varkappa}}\Big),\ \ s\theta\leq 2,
$$ 
there is a constant $C=C(\varkappa,d,r,\theta)>0$ such that
\begin{align}\label{EB0}
\|f\|_{\mL^{s,r}_{t,x}}\leq C\|\nabla f\|^\theta_{\mL^{2,{\varkappa}}_{t,x}}\|f\|_{\mL^{2(1-\theta)s/(2-s\theta),2}_{t,x}}^{1-\theta}.
\end{align}
\el
\begin{proof}
By Gagliado-Nirenberge's interpolation inequality, we have
$$
\|f\|_r\leq C\|\nabla f\|^{\theta}_{\varkappa}\|f\|^{1-\theta}_2.
$$
Since $s\theta\leq 2$, by H\"older's inequality we further have
$$
\|f\|_{\mL^{s,r}_{t,x}}\leq C\|\nabla f\|^{\theta}_{\mL^{2,{\varkappa}}_{t,x}}\|f\|^{1-\theta}_{\mL^{2(1-\theta)s/(2-s\theta),2}_{t,x}}.
$$
The proof is complete.
\end{proof}

Next for fixed $\varkappa\in[1,2]$, we introduce the following index set
$$
\sI_\varkappa:=\Big\{(r,s)\in[2,\infty)\times[1,\infty): \tfrac{1}{2}-\tfrac{1}{r}<\tfrac{1}{s}\big(\tfrac{2}{d}+1-\tfrac{2}{\varkappa}\big)\Big\}.
$$
The following lemma is an easy consequence of \eqref{EB0}.
\bl\label{Le22}
For any $(r,s)\in\sI_\varkappa$ and $\eps\in(0,1)$, there are $\beta\in(1,\infty)$ and constant $C_\eps=C_\eps(r,s,\varkappa,d)>0$ such that
for any $1\leq\tau_1<\tau_2\leq 2$,
\begin{align}\label{EB1}
\|\1_{Q_{\tau_1}}f\|_{\mL^{s,r}_{t,x}}\leq \eps\|\1_{Q_{\tau_2}}\nabla f\|_{\mL^{2,{\varkappa}}_{t,x}}
+C_\eps(\tau_2-\tau_1)^{-1}\|\1_{Q_{\tau_2}}f\|_{\mL^{\beta,2}_{t,x}}.
\end{align}
\el
\begin{proof}
Let $\eta\in C^\infty_c(Q_{\tau_2};[0,1])$ with
$$
\eta|_{Q_{\tau_1}}=1,\ \ |\nabla\eta|\leq 2(\tau_2-\tau_1)^{-1}.
$$
Since $(r,s)\in\sI_\varkappa$, by \eqref{EB0}, there are $\theta\in[0,\frac 2s\wedge 1)$ such that
$$
\|\1_{Q_{\tau_1}}f\|_{\mL^{s,r}_{t,x}}\leq \|\eta f\|_{\mL^{s,r}_{t,x}}\lesssim 
\|\nabla (\eta f)\|^\theta_{\mL^{2,{\varkappa}}_{t,x}}\|\eta f\|_{\mL^{2(1-\theta)s/(2-s\theta),2}_{t,x}}^{1-\theta}.
$$
Moreover, we have
$$
\|\nabla (\eta f)\|_{\mL^{2,{\varkappa}}_{t,x}}
\leq\|\nabla \eta f\|_{\mL^{2,{\varkappa}}_{t,x}}+\|\eta\nabla  f\|_{\mL^{2,{\varkappa}}_{t,x}}
\lesssim(\tau_2-\tau_1)^{-1}\|\1_{Q_{\tau_2}} f\|_{\mL^{2,{\varkappa}}_{t,x}}
+\|\1_{Q_{\tau_2}}\nabla  f\|_{\mL^{2,{\varkappa}}_{t,x}}.
$$
Since $\theta\in[0,1)$ and $s\theta<2$, the desired estimate follows by Young's inequality.
\end{proof}

We need the following elementary variational inequality.
\bl\label{Le21}
Let $N\in\mN$. For any $\alpha_i,p_i\geq 1$ and $\beta_i>0$, $i=1,\cdots,N$, there is a constant $C=C(\alpha_i,\beta_i,p_i,N)>0$ 
such that for all $f_i\in L^\infty([\tau,\sigma])$ and $0<\delta-\tau\leq 1$,
\begin{align}\label{SA3}
\begin{split}
&\inf_{\ell\in C^1([\tau,\delta])}\Bigg\{\sum_{i=1}^N\left(\int^\delta_\tau|\ell'(s)|^{\alpha_i}
|f_i(s)|^{p_i}\dif s\right)^{\frac1{p_i}}\!\!\!:  \ell'(s)\leq 0, \ell(\tau)=1,\ell(\delta)=0\Bigg\}\\
&\qquad\quad\lesssim_C(\delta-\tau)^{-\max\frac{\alpha_i-1}{p_i}-\frac1{\min \beta_i}}\sum_{i=1}^N\left(\int^\delta_\tau |f_i(s)|^{\beta_i}\dif s\right)^{\frac{1}{\beta_i}}.
\end{split}
\end{align}
\el
\begin{proof}
{\it (Step 1).} Fix $\beta\in(0,1)$. We first show that for any $0\leq f\in L^1([\tau,\sigma])$,
\begin{align}\label{SA1}
\begin{split}
\cJ(f)&:=\inf_{\ell\in C^1([\tau,\delta])}\Bigg\{\sum_{i=1}^N\int^\delta_\tau|\ell'(s)|^{\alpha_i}f(s)\dif s: \ell'(s)\leq 0, \ell(\tau)=1,\ell(\delta)=0\Bigg\}\\
&\leq \left(\sum_{i=1}^N2^{\frac{\theta\alpha_i}{\beta}}(\delta-\tau)^{1-\alpha_i-\frac{1}{\beta}}\right)\left(\int^\delta_\tau f(s)^\beta\dif s\right)^{\frac{1}{\beta}},
\end{split}
\end{align}
where $\theta:=1/(\alpha_1\wedge\cdots\alpha_N)$. Let 
$$
\eps:=\left(\frac1{\delta-\tau}\int^\delta_\tau f(s)^\beta\dif s\right)^{\frac{1}{\beta}},\ \ g(s):=f(s)+\eps.
$$
Clearly, we have
$$
\cJ(f)\leq\cJ(g).
$$
Let $g_n:=g*\rho_n$ be the mollifying approximation of $g$. Define
$$
A:=\int^\delta_\tau g(s)^{-\theta}\dif s,\ \  A_n:=\int^\delta_\tau g_n(s)^{-\theta}\dif s,\ \
\ell_n(r):=A^{-1}_n\int^\delta_r g_n(s)^{-\theta}\dif s.
$$
Since $\ell'_n(s)=-A^{-1}_n g_n(s)^{-\theta}\leq 0$, $\ell_n(\tau)=1$ and $\ell_n(\delta)=0$, we have
\begin{align*}
\cJ(g)&\leq \sum_{i=1}^N\int^\delta_\tau|\ell'_n(s)|^{\alpha_i} g(s)\dif s=\sum_{i=1}^NA_n^{-\alpha_i}\int^\delta_\tau g_n(s)^{-\theta\alpha_i}g(s)\dif s.
\end{align*}
Due to $g_n(s)^{-1}\leq \eps^{-1}$ and $1-\theta\alpha_i\leq 0$, by the dominated convergence theorem, we have
\begin{align}
\cJ(g)&\leq \sum_{i=1}^NA^{-\alpha_i}\int^\delta_\tau g(s)^{1-\theta\alpha_i}\dif s
\leq \sum_{i=1}^NA^{-\alpha_i}(\delta-\tau)\eps^{1-\theta\alpha_i}\no\\
&=\sum_{i=1}^NA^{-\alpha_i}(\delta-\tau)^{1-\frac{1-\theta\alpha_i}{\beta}}
\left(\int^\delta_\tau f(s)^\beta\dif s\right)^{\frac{1-\theta\alpha_i}{\beta}}.\label{SA2}
\end{align}
On the other hand, by the inverse H\"older's inequality, we have
\begin{align*}
A^{-1}&=\left(\int^\delta_\tau g(s)^{-\theta}\dif s\right)^{-1}\leq(\delta-\tau)^{-1-\frac{\theta}{\gamma}}\left(\int^\delta_\tau g(s)^\beta\dif s\right)^{\frac{\theta}{\beta}}\\
&\leq (\delta-\tau)^{-1-\frac{\theta}{\beta}}\left(\int^\delta_\tau f(s)^\beta\dif s+(\delta-\tau)\eps^\beta\right)^{\frac{\theta}{\beta}}\\
&= (\delta-\tau)^{-1-\frac{\theta}{\beta}}\left(2\int^\delta_\tau f(s)^\beta\dif s\right)^{\frac{\theta}{\beta}}.
\end{align*}
Substituting this into \eqref{SA2}, we obtain \eqref{SA1}.
\medskip\\
{\it (Step 2).} Let $p:=\max_{i=1,\cdots,N}p_i$. Since $\alpha_i\geq 1$,
by Jensen's inequality with respect to the probability measure $-\1_{[\tau,\delta]}\dif\ell(s)=-\1_{[\tau,\delta]}\ell'(s)\dif s$, we have
$$
\left(\int^\delta_\tau|\ell'(s)|^{\alpha_i}|f_i(s)|^{p_i}\dif s\right)^{\frac{1}{p_i}}
\leq \left(\int^\delta_\tau|\ell'(s)|^{\frac{(\alpha_i-1)p}{p_i}+1}|f_i(s)|^{p}\dif s\right)^{\frac{1}{p}}.
$$
Thus,
\begin{align*}
&\sum_{i=1}^N\left(\int^\delta_\tau|\ell'(s)|^{\alpha_i}
|f_i(s)|^{p_i}\dif s\right)^{\frac1{p_i}}
\leq \sum_{i=1}^N\left(\int^\delta_\tau|\ell'(s)|^{\frac{(\alpha_i-1)p}{p_i}+1}|f_i(s)|^{p}\dif s\right)^{\frac{1}{p}}\\
&\qquad\qquad\leq N\left(\sum_{i=1}^N\int^\delta_\tau|\ell'(s)|^{\frac{(\alpha_i-1)p}{p_i}+1}\left(\sum_{i=1}^N|f_i(s)|^{p}\right)\dif s\right)^{\frac{1}{p}}.
\end{align*}
Let $\sJ$ be the left hand of \eqref{SA3} and $\beta:=\min_{i=1,\cdots,N}\beta_i$. By \eqref{SA1}, we get
\begin{align*}
\sJ&\lesssim_C\left(\sum_{i=1}^N(\delta-\tau)^{-\frac{(\alpha_i-1)p}{p_i}-\frac1\beta}\right)^{\frac1p}
\left(\int^\delta_\tau\left(\sum_{i=1}^N|f_i(s)|^{p}\right)^{\frac\beta{p}}\dif s\right)^{\frac1{\beta}}\\
&\lesssim_C(\delta-\tau)^{-\max\frac{\alpha_i-1}{p_i}-\frac1{p\beta}}\sum_{i=1}^N\left(\int^\delta_\tau|f_i(s)|^{\beta}\dif s\right)^{\frac1{\beta}}\\
&\lesssim_C(\delta-\tau)^{-\max\frac{\alpha_i-1}{p_i}-\frac1{p\beta}}\sum_{i=1}^N(\delta-\tau)^{\frac1\beta-\frac1{\beta_i}}\left(\int^\delta_\tau|f_i(s)|^{\beta_i}\dif s\right)^{\frac1{\beta_i}},
\end{align*}
where the last step is due to $\beta\leq\beta_i$ and H\"older's inequality. The proof is thus complete.
\end{proof}
\br\rm
When $N=1$, the above variational inequality has been used in the proof of \cite[Lemma 2.1]{Be-Sc}. 
For treating the supercritical drifts below, we need a version of $N\geq 2$. The crucial point in \eqref{SA3} 
for us is of course that $\beta_i$ can be smaller than $p_i$ as used in the following 
variational embedding lemma.
\er

The following lemma extends \cite[Lemma 2.1]{Be-Sc} to  time-dependent case and $N\geq 1$.
\bl\label{Le24}
Let $N\in\mN$ and $w=(w_1,\cdots, w_N): \mR\times\mR^d\to\mR^N$ be a bounded measurable function with support in $I\times B_2$, where $I\subset\mR$ is a finite time interval. 
Let  $\alpha_i>0$, $\theta_i\in[0,1]$ and $p_i,q_i,\varkappa_i\geq 1$ satisfy
$$
\alpha_ip_i\geq 1,\ \tfrac{1}{\varkappa_i}=\tfrac{1}{p_i}+\tfrac{\theta_i}{d-1},\ \ i=1,\cdots,N. 
$$
For given $1\leq\tau<\delta\leq 2$ and $Q:=I\times B_\delta$, there are $\gamma,C>0$ depending only on $\alpha_i,\varkappa_i, p_i,d$ such that
\begin{align*}
\cJ(w)&:=\inf_{\eta\in C^1_c(B_\delta;[0,1])}\left\{\sum_{i=1}^N\|w_i|\nabla\eta|^{\alpha_i}\|_{\mL^{p_i,q_i}_{x,t}}: \eta=1\mbox{ on $B_\tau$}\right\}\\
&\lesssim_C(\delta-\tau)^{-\gamma}\sum_{i=1}^N\Big(\|\1_Q\nabla w_i\|^{\theta_i}_{\mL^{\varkappa_i,q_i}_{x,t}}\|\1_Qw_i\|^{1-\theta_i}_{\mL^{\varkappa_i,q_i}_{x,t}}
+\|\1_Qw_i\|_{\mL^{\varkappa_i,q_i}_{x,t}}\Big).
\end{align*}
\el
\begin{proof}
Let 
$$
F_i(x):=\left(\int_I|w_i(t,x)|^{q_i}\dif t\right)^{1/q_i}.
$$ 
For given radial test function $\eta(x)=\ell(|x|)$, by Fubini's theorem and the transform of spherical coordinates, we have
\begin{align*}
\|w_i|\nabla\eta|^{\alpha_i}\|^{p_i}_{\mL^{p_i,q_i}_{x,t}}
&=\int_{\mR^d}F_i^{p_i}|\nabla\eta|^{p_i\alpha_i}
\leq\int^\delta_\tau|\ell'(s)|^{\alpha_i p_i}\left(\int_{\mS^d_s}F_i^{p_i}\right)\dif s,
\end{align*}
where $\mS^d_s:=\{x\in\mR^d:|x|=s\}$.
Since $\frac{1}{\varkappa_i}=\frac{1}{p_i}+\frac{\theta_i}{d-1}$, by the Sobolev embedding in sphere $\mS^d_s$, we have
$$
\|F_i\|_{L^{p_i}(\mS^d_s)}\lesssim \|\nabla F_i\|^{\theta_i}_{L^{\varkappa_i}(\mS^d_s)}\|F_i\|^{1-\theta_i}_{L^{\varkappa_i}(\mS^d_s)}+\|F\|_{L^{\varkappa_i}(\mS^d_s)},\ \ s\in[1,2].
$$
Thus, by Lemma \ref{Le21} with $\beta_i=\varkappa_i$, we have
\begin{align}
\cJ(w)&\lesssim\inf_{\ell\in C^1([\tau,\delta])}\!\!\left\{\sum_{i=1}^N\left(\int^\delta_\tau\!\!|\ell'(s)|^{\alpha_i p_i}\!
\left(\!\int_{\mS^d_s}F^{p_i}\right)\dif s\!\right)^{\frac{1}{p_i}}\!\!\!\!: \ell'(s)\leq 0, \ell(\tau)=1,\ell(\delta)=0\right\}\no\\
&\lesssim (\delta-\tau)^{-\gamma}\sum_{i=1}^N\left(\int^\delta_\tau \left(\int_{\mS^d_s}F^{p_i}\right)^{\frac{\varkappa_i}{p_i}}\dif s\right)^{\frac{1}{\varkappa_i}}\no\\
&\lesssim (\delta-\tau)^{-\gamma}\sum_{i=1}^N\left(\int^\delta_\tau 
\left(\|\nabla F_i\|^{\theta_i}_{L^{\varkappa_i}(\mS^d_s)}\|F\|^{1-\theta_i}_{L^{\varkappa_i}(\mS^d_s)}
+\|F\|_{L^{\varkappa_i}(\mS^d_s)}\right)^{\varkappa_i}\dif s\right)^{\frac{1}{\varkappa_i}}\no\\
&\lesssim(\delta-\tau)^{-\gamma}\sum_{i=1}^N
\left(\|\nabla F_i\|^{\theta_i}_{L^{\varkappa_i}(B_\delta)}\|F_i\|^{1-\theta_i}_{L^{\varkappa_i}(B_\delta)}
+\|F_i\|_{L^{\varkappa_i}(B_\delta)}\right).\label{KP1}
\end{align}
On the other hand, let 
$$
F^{(\eps)}_i(x):=\left(\int_I(|w_i(t,x)|^{q_i}+\eps)\dif t\right)^{1/q_i}.
$$ 
By the chain rule and H\"older's inequality, we have
\begin{align*}
|\nabla F^{(\eps)}_i(x)|&\leq\left(\int_I(|w_i(t,x)|^{q_i}+\eps)\dif t\right)^{\frac{1-q_i}{q_i}}\!\!\!\!\int_I |w_i(t,x)|^{q_i-1}|\nabla w_i(t,x)|\dif t\\
&\leq \left(\int_I(|w_i(t,x)|^{q_i}+\eps)\dif t\right)^{\frac{1-q_i}{q_i}}\!\!\!\!
\left(\int_I|w_i(t,x)|^{q_i}\dif t\right)^{\frac{q_i-1}{q_i}}\!\!\!\!\|\nabla w_i(\cdot,x)\|_{L^{q_i}(I)}.
\end{align*}
Letting $\eps\downarrow 0$, we obtain
$$
|\nabla F_i(x)|\leq \|\nabla w_i(\cdot,x)\|_{L^{q_i}(I)}.
$$
Substituting this into \eqref{KP1}, we obtain the desired estimate. 
\end{proof}
\br\rm
Suppose that $w:\mR^d\to\mR$ is time-independent and $\frac{1}{\varkappa}\leq \frac{1}{p}+\frac{1}{d}$.  
For $N=1$, by using Sobolev's embedding directly, we have
\begin{align*}
\cJ(w)\leq (\delta-\tau)^{-\alpha}\|w\|_{L^p(B_\delta)}
\lesssim_C (\delta-\tau)^{-\alpha}\Big(\|\nabla w\|_{L^\varkappa(B_\delta)}+\|w\|_{L^\varkappa(B_\delta)}\Big).
\end{align*}
However, by Lemma \ref{Le22}, we have for $\frac{1}{\varkappa}\leq\frac{1}{p}+\frac{1}{d-1}$,
\begin{align*}
\cJ(w)\lesssim_C (\delta-\tau)^{-\gamma}\Big(\|\nabla w\|_{L^\varkappa(B_\delta)}+\|w\|_{L^\varkappa(B_\delta)}\Big),
\end{align*}
which clearly has smaller $\varkappa$ than the above estimate. 
\er

We also need the following iteration lemma (cf. \cite[Lemma 4.3]{Ha-Li}).
\bl\label{Le26}
Let $h(\tau)\geq 0$ be bounded in $[\tau_1,\tau_2]$ with $\tau_1\geq 0$. Let $A,B>0$. Suppose that for some $\alpha\geq 0$, $\theta\in(0,1)$ and any
$\tau_1\leq\tau<\tau'\leq\tau_2$, 
$$
h(\tau)\leq\theta h(\tau')+(\tau'-\tau)^{-\alpha}A+B.
$$
Then there is a constant $C=C(\alpha,\theta)>0$ such that
$$
h(\tau_1)\lesssim_C((\tau_2-\tau_1)^{-\alpha}A+B).
$$
\el
\section{Maximum principle for linear parabolic equations}

Let $p_0\in(\frac{d}{2},\infty]$ and $p_1\in[1,\infty]$ be as in \eqref{PP0}. We define $\varkappa\in[1,2]$ by
\begin{align}\label{DW2}
\tfrac{2}{{\varkappa}}=\tfrac{1}{{p_0}}+1.
\end{align}
For a set $Q\subset\mR^{d+1}$, we also introduce
$$
\sV_Q:=\Big\{f\in \mL^1_{loc}: \|f\|_{\sV_Q}:=\|\1_Q f\|_{\mL^{\infty,2}_{t,x}}+\|\1_Q\nabla_xf\|_{\mL^{{\varkappa},2}_{x,t}}<\infty\Big\}
$$
and 
\begin{align}\label{VV}
\widetilde\sV:=\Big\{f\in \mL^1_{loc}: 
\nor f\nor_{\widetilde\sV}:=\nor f\nor_{\widetilde\mL^{\infty,2}_{t,x}}+\nor\nabla_xf\nor_{\widetilde\mL^{{\varkappa},2}_{x,t}}<\infty\Big\}.
\end{align}
\subsection{Energy type estimate}
In this subsection we fix $1\leq\tau_1<\tau_2\leq 2$ and 
$$
Q_i:=Q_{\tau_i}=[-\tau^2_i,\tau^2_i]\times B_{\tau_i},\ i=1,2.
$$
Let $\sC$ be the class of all functions with the following form:
$$
\eta(t,x)=\eta_{\rm t}(t)\eta_{\rm x}(x),
$$
where $\eta_{\rm x}\in C^1_c(B_{\tau_2};[0,1])$ and $\eta_{\rm t}\in C^1_c([-\tau^2_2,\tau^2_2];[0,1])$  with 
\begin{align}\label{Eta}
\eta_{\rm x}|_{B_{\tau_1}}=1,\  \eta_{\rm t}|_{[-\tau_1^2,\tau^2_1]}=1,\  \|\p_t\eta_{\rm t}\|_\infty\leq 4(\tau_2-\tau_1)^{-1}.
\end{align}
\newpage
We first prepare the following important variational estimate.
\bl\label{Le35}
Let $g_1\in\widetilde\mL^{p_1,\infty}_{x,t}$ and $g_2\in\widetilde\mL^{p_3,q_3}_{x,t}$, where $1\leq p_3\leq q_3\leq\infty$ satisfy \eqref{Re01}.
For any $\eps\in(0,1)$, there are $(r,s)\in\sI_\varkappa$ and constants $\gamma>0$ and $C_\eps$ depending on the norm of $\|\1_{Q_2}g_1\|_{\mL^{p_1,\infty}_{x,t}}$
and $\|\1_{Q_2}g_2\|_{\mL^{p_3,q_3}_{x,t}}$ such that for any $w\in\sV_{Q_2}$,
\begin{align}\label{ES0}
\begin{split}
&\inf_{\eta\in\sC}\Big(\|g_1 w^2|\nabla\eta|^2\|_{\mL^{1,1}_{t,x}}
+\|g_2 w^2|\nabla\eta|\|_{\mL^{1,1}_{t,x}}\Big)\\
&\qquad\leq \eps\|w\|^2_{\sV_{Q_2}}+C_\eps(\tau_2-\tau_1)^{-\gamma}\|\1_{Q_2}w\|_{\mL^{s,r}_{t,x}}^2.
\end{split}
\end{align}
\el
\begin{proof}
Let $\bar p_1=\frac{p_1}{p_1-1}$. By H\"older's inequality, we have
$$
\|g_1 w^2|\nabla\eta|^2\|_{\mL^{1,1}_{t,x}}
\leq\|\1_{Q_2}g_1\|_{\mL^{p,\infty}_{x,t}}\|w^{2}|\nabla\eta|^2\|_{\mL^{\bar p,1}_{x,t}}
\leq\|\1_{Q_2}g_1\|_{\mL^{p_1,\infty}_{x,t}}\|w|\nabla\eta|\|^2_{\mL^{2\bar p_1,2}_{x,t}}.
$$
By \eqref{Re01}, one can choose $\delta\in(0,1)$ so that
$$
\tfrac{d-1}2\big(\tfrac1{p_3}-\tfrac1{q_3}\big)\vartheta_1<\delta<1-\tfrac{\vartheta_1+2+d(\vartheta_2-\vartheta_1)}{2q_3}.
$$
 Let  $r,s\in(1,\infty]$  be defined by
$$
\tfrac1{q_3}+\tfrac{1}{s}+\delta=1,\ \  \tfrac1{p_3}+\tfrac{1}{r}+\tfrac{1}{s}=1,
$$
and let
$$
s=r:=2(1-\delta)s,\ \ \bar p_3:=\delta r.
$$
By H\"older's inequality, we have
\begin{align*}
\|g_2 w^2\nabla\eta\|_{\mL^{1,1}_{t,x}}
&\leq\|\1_{Q_2}g_2\|_{\mL^{p_3,q_3}_{x,t}}\|\1_{Q_2}|w|^{2(1-\delta)}\|_{\mL^{s,s}_{x,t}}\||w|^{2\delta}\nabla\eta\|_{\mL^{r,1/\delta}_{x,t}}\no\\
&=\|\1_{Q_2}g_2\|_{\mL^{p_3,q_3}_{x,t}}\|\1_{Q_2}w\|^{2(1-\delta)}_{\mL^{s,r}_{t,x}}\|w|\nabla\eta|^{\frac{1}{2\delta}}\|^{2\delta}_{\mL^{2\bar p_3,2}_{x,t}}\no\\
&\leq\|\1_{Q_2}g\|^{1/(1-\delta)}_{\mL^{p_3,q_3}_{x,t}}\|\1_{Q_2}w\|^{2}_{\mL^{s,r}_{t,x}}
+\|w|\nabla\eta|^{\frac{1}{2\delta}}\|^{2}_{\mL^{2\bar p_3,2}_{x,t}}.
\end{align*}
Thus, if we denote by $\cJ(w)$ the left hand of \eqref{ES0}, then
\begin{align*}
\cJ(w)\lesssim \|\1_{Q_2}w\|^{2}_{\mL^{s,r}_{t,x}}+\inf_{\eta\in\sC}\Big(\|w|\nabla\eta|\|^2_{\mL^{2\bar p_1,2}_{x,t}}
+\|w|\nabla\eta|^{\frac{1}{2\delta}}\|^{2}_{\mL^{2\bar p_3,2}_{x,t}}\Big).
\end{align*}
Note that
$$
\delta<1-\tfrac{\vartheta_1+2+d(\vartheta_2-\vartheta_1)}{2q_3}\Rightarrow (r,s)\in\sI_\varkappa
$$
and
$$
\tfrac{d-1}2\big(\tfrac1{p_3}-\tfrac1{q_3}\big)\vartheta_1<\delta\Rightarrow \tfrac{1}{2p_0}+\tfrac12=\tfrac{1}{\varkappa}<\tfrac{1}{2\bar p_3}+\tfrac{1}{d-1}.
$$
One can choose $\theta_1,\theta_3\in(0,1)$ being close to $1$ so that
$$
\tfrac{1}{\varkappa}=\tfrac{1}{2\bar p_i}+\tfrac{\theta_i}{d-1},\ \ i=1,3.
$$
Thus by  Lemma \ref{Le24} and Young's inequality, we have for any $\eps>0$ and some $\gamma>0$,
\begin{align*}
\cJ(w)&\leq C\|\1_{Q_2}w\|^{2}_{\mL^{s,r}_{t,x}}+
\eps\|\1_{Q_2}\nabla w\|^2_{\mL^{\varkappa,2}_{x,t}}
+C_\eps (\tau_2-\tau_1)^{-\gamma}\|\1_{Q_2}w\|^2_{\mL^{\varkappa,2}_{x,t}}.
\end{align*}
Note that
$$
\|\1_{Q_2}w\|_{\mL^{\varkappa,2}_{x,t}}\lesssim\|\1_{Q_2}w\|_{\mL^{2,2}_{t,x}}\lesssim\|\1_{Q_2}w\|_{\mL^{s,r}_{t,x}}.
$$
The proof is thus complete.
\end{proof}

Recall $\Theta$ being the parameter set \eqref{TH}.
Now we prove the following local energy estimate by Lemma \ref{Le35}.
\bl\label{Le32}
Under {\bf (H$^a$)} and {\bf (H$^b$)}, for any 
$f\in \mL^{q_4,p_4}_{t,x}(Q_2)$ with $(p_4,q_4)\in\mI^d_{p_0}$, there are $(r_i,s_i)\in\sI_\varkappa$, $i=1,2,3$, $\gamma=\gamma(p_0,p_1,d)\geq1$ and constant $C=C(\Theta)>0$
such that for any Lipschitz weak subsolution $u$ of PDE \eqref{PDE0} and $t\geq 0$,
$$
\|w\cI^t\|^2_{\sV_{Q_1}}\lesssim_C(\tau_2-\tau_1)^{-\gamma}\sum_{i=1,2}\|\1_{Q_2}w\cI^t\|^2_{\mL^{s_i,r_i}_{t,x}}
+\|f\1_{Q_2}\|^2_{\mL^{q_4,p_4}_{t,x}}\|\1_{\{w\not=0\}\cap Q_2}\cI^t\|^2_{\mL^{s_3,r_3}_{t,x}},
$$
where
$\cI^t(\cdot):=\1_{(-\infty,t]}(\cdot)$ and $w:=(u-\kappa)^+$ and $\kappa\geq 0.$ 
\el
\begin{proof}
We divide the proof into three steps.
\medskip\\
(i) Fix $\eta\in\sC$ (see \eqref{Eta}). In this step we show that for all $t\in\mR$,
\begin{align}\label{HA1}
\begin{split}
\|(\eta w)(t)\|^2_2&\leq\<\!\<\p_s\eta^2, w^2\cI^t\>\!\>- 2\<\!\<a\cdot\nabla u,\nabla(\eta^2w)\cI^t\>\!\>\\
&\quad+ 2\<\!\<b\cdot \nabla u, \eta^2w\cI^t\>\!\>+ 2\<\!\<f,\eta^2w\cI^t\>\!\>.
\end{split}
\end{align}
Since we want to take the test function $\varphi=w\eta^2$ in \eqref{Def0}, and $\p_s u$ only makes sense in the distributional sense, 
we shall first approximate $u$ by its Steklov mean:
\begin{align}\label{Def2}
S_h u(t,x):=\frac{1}{h}\int_0^{h} u(t+s,x) \dif s=\frac{1}{h}\int_t^{t+h} u(s,x) \dif s,\ \ h\in(0,1).
\end{align}
Let $u_h:=S_h u$ and $S^*_h$ be the adjoint operator of $S_h$. Let $\varphi$ be a nonnegative Lipschitz function in $\mR^{d+1}$ with compact support in $Q_2$.
By Definition \ref{Def1} with test function $S^*_h\varphi$, using integration by parts and Fubini's theorem,
one sees that
\begin{align}\label{Def11}
\<\!\<\p_s u_h,\varphi\>\!\>\leq-\<\!\<S_h(a\cdot\nabla u),\nabla\varphi\>\!\>+\<\!\<S_h(b\cdot\nabla u),\varphi\>\!\>+\<\!\<f_h,\varphi\>\!\>.
\end{align}
Now fix $t\in\mR$ and define
\begin{align*}
\psi_t^{\eps}(s)=\1_{(-\infty,t]}(s)+(1-\eps^{-1}(s-t))\1_{(t,t+\eps]}(s), \ \eps\in(0,1).
\end{align*}
Let $w_h:= (u_h-k)^+$. Note that
\begin{align*}
&2\<\p_s u_h, w_h\eta^2 \psi_t^{\eps}\>= 2\<\p_s w_h, w_h\eta^2 \psi_t^{\eps}\>\\
&\quad=\p_s\<w_h^2,\eta^2\psi_t^\eps\>-\int_{\mR^d}  w_h^2 \eta^2 (\psi_t^\eps)'-\int_{\mR^d} w_h^2 (\p_s \eta^2 \psi_t^{\eps}).
\end{align*}
By \eqref{Def11} with $\varphi = w_h \eta^2 \psi_{t, \eps}$ and $\int_{\mR^{d+1}}\p_s\<w_h^2,\eta^2\psi_t^\eps\>=0$, we get
\begin{align*}
-\int_{\mR^{d+1}}\eta^2 w^2_h(\psi_t^{\eps})' \leq&  \int_{\mR^{d+1}} w^2_h(\p_s \eta^2\psi_t^{\eps}) -
2\<\!\<S_h(a\cdot \nabla u), \nabla (w_h \eta^2\psi_t^{\eps})\>\!\> \\
& +2\<\!\<S_h(b\cdot \nabla u), w_h \eta^2\psi_t^{\eps}\>\!\> + 2\<\!\<f_h, w_h \eta^2\psi_t^{\eps}\>\!\>.
\end{align*}
Letting  $h\downarrow 0$ and by the dominated convergence theorem, we obtain 
\begin{align*}
-\int_{\mR^{d+1}}(\eta w)^2 (\psi_t^{\eps})' \leq&  \int_{\mR^{d+1}} w^2(\p_t \eta^2\psi_t^{\eps}) -
2\<\!\<a\cdot \nabla u, \nabla (w \eta^2\psi_t^{\eps})\>\!\> \\
& +2\<\!\<b\cdot \nabla u, w\eta^2\psi_t^{\eps}\>\!\> + 2\<\!\<f, w \eta^2\psi_t^{\eps}\>\!\>.
\end{align*}
Since $\lim_{\eps\downarrow 0}\psi_t^{\eps}(s)=\cI^t(s)$ for each $s\in\mR$,
the right hand side of the above inequality converges to the right hand side of 
\eqref{HA1} as $\eps\downarrow 0$. 
On the other hand, by the Lebesgue differential theorem, we also have
$$
-\int_{\mR^{d+1}}(\eta w)^2(\psi_t^{\eps})' =\frac{1}{\eps}\int_t^{t+\eps}\|(\eta w)(s)\|_2^2 \dif s\stackrel{\eps\downarrow 0}{ \to}\|(\eta w)(t)\|_2^2.
$$
Thus, we obtain \eqref{HA1}. 
\medskip\\
(ii) Recalling $w=(u-\kappa)^+$ and noting that
\begin{align}\label{HA2}
\nabla u\cdot\nabla w=|\nabla w|^2,\ (\nabla u) w=(\nabla w)w=\nabla w^2/2,
\end{align}
by the chain rule and Young's inequality, we have
\begin{align*}
-\<a\cdot\nabla u,\nabla(\eta^2w)\>
&=-\int_{\mR^d}\eta^2(\nabla w)^*a\nabla w-2\int_{\mR^d} \eta w (\nabla\eta)^*a\nabla w\\
&\leq-\frac{1}{2}\int_{\mR^d}\eta^2(\nabla w)^*a\nabla w+4\int_{\mR^d} w^2|\nabla\eta|^2\frac{|a\nabla w|^2}{(\nabla w)^*a\nabla w}\\
&\stackrel{\eqref{DW3}}{\leq} -\frac{1}{2}\int_{\mR^d}\eta^2|\nabla w|^2\lambda +4\int_{\mR^d} w^2|\nabla\eta|^2\mu,
\end{align*}
which in turn gives that
\begin{align}\label{HQ1}
-\<\!\<a\cdot\nabla u,\nabla(\eta^2w)\cI^t\>\!\>
\leq -\tfrac{1}{2}\|\eta\nabla w\lambda^{\frac12}\cI^t\|^2_{\mL^{2,2}_{x,t}}+4\|\mu|\nabla\eta|^2w^2\cI^t\|_{\mL^{1,1}_{t,x}}.
\end{align}
Due to $b=b_1+b_2$ and $(\div b_2)^-=0$, by \eqref{HA2} and the integration by parts, we have
\begin{align*}
\<\!\<b\cdot\nabla u, \eta^2w\cI^t\>\!\>
&=\<\!\<\eta b_1\cdot\nabla w\lambda^{\frac12}, \lambda^{-\frac12}\eta w\cI^t\>\!\>+\tfrac{1}{2}\<\!\<b_2\cdot\nabla w^2, \eta^2\cI^t\>\!\>\\
&\leq
\tfrac{1}{4}\|\eta \nabla w\lambda^{\frac12}\cI^t\|^2_{\mL^{2,2}_{t,x}}
+4\|\lambda^{-\frac12} b_1\eta w\cI^t\|^2_{\mL^{2,2}_{t,x}}-\<\!\<b_2\cdot \nabla \eta, w^2\cI^t\>\!\>.
\end{align*}
Let $(r_2,s_2)$ be defined by 
$$
\tfrac{1}{2p_0}+\tfrac{1}{p_2}+\tfrac{1}{r_2}=\tfrac{1}{2},\ \ \tfrac{1}{q_2}+\tfrac{1}{s_2}=\tfrac{1}{2},
$$
which belongs to $\sI_\varkappa$ by  \eqref{Re1}.
By H\"older's inequality, we have
$$
\|\lambda^{-\frac12} b_1\eta w\cI^t\|^2_{\mL^{2,2}_{t,x}}\leq
\|\lambda^{-1}\1_{B_{\tau_2}}\|_{p_0}\|b_1\1_{Q_2} \|^2_{\mL^{q_2,p_2}_{t,x}}\|\eta w\cI^t\|^2_{\mL^{s_2,r_2}_{t,x}}.
$$
Hence, 
\begin{align}\label{HQ2}
\<\!\<b\cdot\nabla u, \eta^2w\cI^t\>\!\>
&\leq\tfrac{1}{4}\|\eta \nabla w\lambda^{\frac12}\cI^t\|^2_{\mL^{2,2}_{t,x}}
+C\|\eta w\cI^t\|^2_{\mL^{s_2,r_2}_{t,x}}+\|b_2\cdot\nabla \eta w^2\cI^t\|_{\mL^{1,1}_{t,x}}.
\end{align}
Similarly, let $(r_3,s_3)$  be defined by
$$
\tfrac{1}{p_4}+\tfrac{1}{r_3}=\tfrac{1}{2},\ \ \tfrac{1}{q_4}+\tfrac{1}{s_3}=1,
$$
which belongs to $\sI_\varkappa$ by $(p_4,q_4)\in\mI^d_{p_0}$.
By H\"older's inequality and Young's inequality again, we also have
\begin{align}\label{HQ3}
\<\!\<f, \eta^2w\cI^t\>\!\>
&\leq\|f\eta\|_{\mL^{q_4,p_4}_{t,x}}\|\eta w\cI^t\|_{\mL^{\infty,2}_{t,x}}\|\1_{\{w\not=0\}\cap Q_2}\cI^t\|_{\mL^{s_3,r_3}_{t,x}}\no\\
&\leq \tfrac{1}{4}\|\eta w\cI^t\|_{\mL^{\infty,2}_{t,x}}^2+4\|f\eta\|^2_{\mL^{q_4,p_4}_{t,x}}\|\1_{\{w\not=0\}\cap Q_2}\cI^t\|^2_{\mL^{s_3,r_3}_{t,x}}.
\end{align}
Combining \eqref{HA1} and \eqref{HQ1}-\eqref{HQ3}, we obtain
\begin{align}
\|(\eta w)(t)\|^2_2&\leq 2\|\p_s\eta w^2\cI^t\|_{\mL^{1,1}_{t,x}}
+ 8\|\mu|\nabla\eta|^2w^2\cI^t\|_{\mL^{1,1}_{t,x}}-\tfrac{1}{2}\|\eta \nabla w\lambda^{\frac12}\cI^t\|^2_{\mL^{2,2}_{t,x}}\no\\
&\quad+2\|b_2\cdot\nabla \eta w^2\cI^t\|_{\mL^{1,1}_{t,x}}+C\|\eta w\cI^t\|^2_{\mL^{s_2,r_2}_{t,x}}\no\\
&\quad+\tfrac{1}{2}\|\eta w\cI^t\|_{\mL^{\infty,2}_{t,x}}^2
+8\|f\eta\|^2_{\mL^{q_4,p_4}_{t,x}}\|\1_{\{w\not=0\}\cap Q_2}\cI^t\|^2_{\mL^{s_3,r_3}_{t,x}}.\label{DC0}
\end{align}
(iii) By \eqref{DW2} and H\"older's inequality, we have
$$
\|\eta \nabla w\cI^t\|^2_{\mL^{\chi,2}_{t,x}}\leq \nor\lambda^{-1}\nor_{p_0}\|\eta \nabla w\lambda^{\frac12}\cI^t\|^2_{\mL^{2,2}_{t,x}}.
$$
Since $\eta|_{Q_1}=1$ and $\eta|_{Q_2^c}=0$, by  \eqref{DC0} and \eqref{Eta}, we obtain that for any $\eta\in\sC$,
\begin{align}\label{DD4}
\begin{split}
\|w\cI^t\|^2_{\sV_{Q_1}}&\leq
\|\eta w\cI^t\|_{\mL^{\infty,2}_{t,x}}^2+\|\eta \nabla w\cI^t\|^2_{\mL^{\chi,2}_{t,x}}\\
&\lesssim_C(\tau_2-\tau_1)^{-1}\|\1_{Q_2} w\cI^t\|^2_{\mL^{2,2}_{t,x}}+\cJ(w,\eta)+\|\1_{Q_2} w\cI^t\|^2_{\mL^{s_2,r_2}_{t,x}}\\
&\quad+\|f\1_{Q_2}\|^2_{\mL^{q_4,p_4}_{t,x}}\|\1_{\{w\not=0\}\cap Q_2}\cI^t\|^2_{\mL^{s_3,r_3}_{t,x}},
\end{split}
\end{align}
where
$$
\cJ(w,\eta):=\|\mu|\nabla\eta|^2w^2\cI^t\|_{\mL^{1,1}_{t,x}}+\|b_2\cdot\nabla\eta w^2\cI^t\|_{\mL^{1,1}_{t,x}}.
$$
By Lemma \ref{Le35}, there are $(r_1,s_1)\in\sI_\varkappa$ and $\gamma_0>0$ such that for all $\eps\in(0,1)$,
$$
\inf_{\eta\in\sC}\cJ(w,\eta)\leq\eps\|w\cI^t\|^2_{\sV_{Q_2}}+C_\eps
(\tau_2-\tau_1)^{-\gamma_0}\|\1_{Q_2}w\cI^t\|^2_{\mL^{s_1,r_1}_{t,x}}.
$$
Substituting this into \eqref{DD4}  and by 
$\|\1_{Q_2}w\|_{\mL^{2,2}_{t,x}}\lesssim\|\1_{Q_2}w\|_{\mL^{s_1,r_1}_{t,x}},$ 
we get that for some $\gamma=\gamma(p_0,p_1,d)\geq 1$, $C>0$ and all $1\leq\tau_1<\tau_2\leq 2$,
\begin{align*}
\|w\cI^t\|^2_{\sV_{Q_1}}&\leq
\tfrac12\|w\cI^t\|^2_{\sV_{Q_2}}+C(\tau_2-\tau_1)^{-\gamma}\sum_{i=1,2}\|\1_{Q_2} w\cI^t\|^2_{\mL^{s_i,r_i}_{t,x}}\\
&\quad+\|f\1_{Q_2}\|^2_{\mL^{q_4,p_4}_{t,x}}\|\1_{\{w\not=0\}\cap Q_2}\cI^t\|^2_{\mL^{s_3,r_3}_{t,x}}.
\end{align*}
Recall $Q_i=Q_{\tau_i}$ for $i=1,2$. If we let $h(\tau):=\|w\cI^t\|^2_{\sV_{Q_\tau}}$, then the above inequality implies that
for any $\tau_1\leq \tau<\tau'\leq\tau_2$,
\begin{align*}
h(\tau)&\leq
\tfrac12 h(\tau')+C(\tau'-\tau)^{-\gamma}\sum_{i=1,2}\|\1_{Q_{\tau_2}} w\cI^t\|^2_{\mL^{s_i,r_i}_{t,x}}\\
&\quad+\|f\1_{Q_{\tau_2}}\|^2_{\mL^{q_4,p_4}_{t,x}}\|\1_{\{w\not=0\}\cap Q_{\tau_2}}\cI^t\|^2_{\mL^{s_3,r_3}_{t,x}}.
\end{align*}
The desired estimate now follows by Lemma \ref{Le26}.
\end{proof}

\subsection{Local maximum estimate}
The following lemma is easy by H\"older's inequality.
\bl\label{Le12}
Let $Q=I\times D\subset \mR\times\mR^{d}$ 
be a bounded domain. For any $p,q\in[1,\infty)$, there are constants $C_1,C_2>0$ only depending on $Q, p,q$ such that for any $A\subset Q$,
$$
\|\1_A\|_{\mL^{q,p}_{t,x}}+\|\1_A\|_{\mL^{p,q}_{x,t}}\leq C_1\|\1_A\|^{1/(p\vee q)}_{\mL^{1,1}_{t,x}}
\leq C_2\big(\|\1_A\|_{\mL^{q,p}_{t,x}}+\|\1_A\|_{\mL^{p,q}_{x,t}}\big)^{1/(p\vee q)}.
$$
\el

We need the following simple De-Giorgi's iteration lemma.
\bl\label{Le15}
Let $(a_n)_{n\in\mN}$ be a sequence of nonnegative real numbers. Suppose that for some $C_0,\lambda>1$ and $\delta_j>0$, $j=1,\cdots,m$,
$$
a_{n+1}\leq C_0\lambda^n a_n\sum_{j=1}^m a_n^{\delta_j},\ \ n=1,2,\cdots.
$$
If $a_1\leq (mC_0\lambda^{(1+\delta)/\delta})^{-1/\delta}$, where $\delta=\delta_1\wedge\cdots\wedge\delta_m$, then
$$
\lim_{n\to\infty} a_n=0.
$$
\el
\begin{proof}
We use induction to prove that if $a_1\leq  (mC_0\lambda^{(1+\delta)/\delta})^{-1/\delta}\leq 1$, then
$$
a_n\leq a_1\lambda^{-(n-1)/\delta},\ \ \forall n\in\mN.
$$
By the induction hypothesis, we have
\begin{align*}
a_{n+1}&\leq mC_0\lambda^n a^{1+\delta}_n\leq mC_0 a_1^{1+\delta}\lambda^{n-(n-1)(1+\delta)/\delta}\\
&=\big(mC_0 a_1^\delta\lambda^{(1+\delta)/\delta}\big)a_1\lambda^{-n/\delta}\leq a_1\lambda^{-n/\delta},
\end{align*}
where the last step is due to $mC_0 a_1^\delta\lambda^{(1+\delta)/\delta}\leq 1$.
\end{proof}
\bl
Let $1\leq \tau_1<\tau_0\leq 2$ and $0<\kappa_0<\kappa_1$. Define
$$
\Gamma_i:=Q_{\tau_i},\ \ w_i:=(u-\kappa_i)^+,\ \ i=0,1.
$$
(i) For any $r,s\in[1,\infty]$, we have
\begin{align}\label{LK1}
\|\1_{\{w_1\not=0\}\cap\Gamma_0}\|_{\mL^{s,r}_{t,x}}\leq \|w_0\1_{\Gamma_0}\|_{\mL^{s,r}_{t,x}}/(\kappa_1-\kappa_0).
\end{align}
(ii) For any $r\in[1,{\varkappa})$ and $s\in[1,2)$, there is a universal constant $C>0$ such that
\begin{align}\label{LK12}
\|\1_{\Gamma_0} \nabla w_1\|_{\mL^{r,s}_{x,t}}
\leq C\|w_1\|_{\sV_{\Gamma_0}}\Big(\tfrac{\|w_0\1_{\Gamma_0}\|_{\mL^{1,1}_{t,x}}}{\kappa_1-\kappa_0}\Big)
^{(\frac1 r-\frac1\varkappa)\wedge(\frac1s-\frac12)}.
\end{align}
(iii) For any $(r,s)\in\sI_\varkappa$, 
there are $\delta\in(0,1)$ and $C=C(r,s,d,\varkappa)>0$  such that
\begin{align}\label{LK4}
\|\1_{\Gamma_1}w_1\|_{\mL^{s,r}_{t,x}}\leq C (\tau_0-\tau_1)^{-1}\|w_1\|_{\sV_{\Gamma_0}}
\cdot\Big(\tfrac{\|w_0\1_{\Gamma_0}\|_{\mL^{1,1}_{t,x}}}{\kappa_1-\kappa_0}\Big)^{\delta}.
\end{align}
\el
\begin{proof}
(i) Noting that
$$
w_0|_{w_1\not=0}=(u-\kappa_1+\kappa_1-\kappa_0)^+|_{w_1\not=0}\geq \kappa_1-\kappa_0,
$$
for given $r,s\in[1,\infty]$, we have
$$
\|w_0\1_{\Gamma_0}\|_{\mL^{s,r}_{t,x}}
\geq\|w_0\1_{\{{w_1\not=0}\}\cap\Gamma_0}\|_{\mL^{s,r}_{t,x}}\geq (\kappa_1-\kappa_0) \|\1_{\{w_1\not=0\}\cap\Gamma_0}\|_{\mL^{s,r}_{t,x}},
$$
which implies \eqref{LK1}.

(ii) Let $\frac{1}{r}=\frac{1}{{\varkappa}}+\frac{1}{r'}$
and $\frac{1}{s}=\frac{1}{2}+\frac{1}{s'}$. By H\"older's inequality, we have
\begin{align*}
\|\1_{\Gamma_0} \nabla w_1\|_{\mL^{r,s}_{x,t}}
&=\|\1_{\Gamma_0\cap\{w_1\not=0\}} \nabla w_1\|_{\mL^{r,s}_{x,t}}
\leq\|\1_{\Gamma_0} \nabla w_1\|_{\mL^{{\varkappa},2}_{x,t}}\|\1_{\Gamma_0\cap\{w_1\not=0\}}\|_{\mL^{r',s'}_{x,t}}\\
&\lesssim\|w_1\|_{\sV_{\Gamma_0}}\|\1_{\Gamma_0\cap\{w_1\not=0\}}\|^{1/(s'\vee r')}_{\mL^{1,1}_{t,x}},
\end{align*}
which implies \eqref{LK12} by \eqref{LK1}.

(iii) Since $\tfrac{1}{2}-\tfrac{1}{r}<\tfrac{1}{s}(\tfrac{2}{d}-\tfrac{1}{p_0})$,  we can choose $r',\beta>r$ and $s',\theta>s$ so that
$$
\tfrac{1}{r'}+\tfrac{1}{\beta}=\tfrac{1}{r},\quad \tfrac{1}{s'}+\tfrac{1}{\theta}=\tfrac{1}{s},\ \ 
\tfrac{1}{2}-\tfrac{1}{r'}<\tfrac{1}{s'}(\tfrac{2}{d}-\tfrac{1}{p_0}).
$$
By  H\"older's inequality, we have
$$
\|w_1\1_{\Gamma_1}\|_{\mL^{s,r}_{t,x}}\leq\|w_1\1_{\Gamma_1}\|_{\mL^{s',r'}_{t,x}}\|\1_{\{w_1\not=0\}\cap\Gamma_1}\|_{\mL^{\theta,\beta}_{t,x}},
$$
and by Lemma \ref{Le22} and \eqref{HQ5},
\begin{align*}
\|w_1\1_{\Gamma_1}\|_{\mL^{s',r'}_{t,x}}&\lesssim_C\|\1_{\Gamma_0}\nabla w_1\|_{\mL^{2,{\varkappa}}_{t,x}}
+(\tau_0-\tau_1)^{-1}\|\1_{\Gamma_0}w_1\|_{\mL^{\infty,2}_{t,x}}\\
&\lesssim_C(\tau_0-\tau_1)^{-1}\|w_1\|_{\sV_{\Gamma_0}},
\end{align*}
which in turn yields \eqref{LK4} by Lemma \ref{Le12} and \eqref{LK1}.
\end{proof}

Now we can show the following local maximum principle for PDE \eqref{PDE0}.
\bt[Local maximum estimate]
Under the assumption of Lemma \ref{Le32}, for any $p>0$, there is a constant $C=C(p,\Theta)>0$ such that 
for any Lipschitz weak subsolution $u$ of PDE \eqref{PDE0},
\begin{align}\label{JK2}
\|u^+\1_{Q_1}\|_{\infty}\leq C\left(\|u^+\1_{Q_2}\|_{\mL^{p,p}_{t,x}}+\|f\1_{Q_2}\|_{\mL^{q_4,p_4}_{t,x}}\right),
\end{align}
where $Q_1:=[-1,1]\times b_1$ and $Q_2:=[-4,4]\times b_2$.
\et
\begin{proof}
Fix $1\leq\tau<\sigma\leq 2$.  Let $\kappa>0$, which will be determined below. For $n\in\mN$, define
$$
\tau_n=\tau+(\sigma-\tau)2^{1-n},\ \ \tilde\tau_n:=\tau+3(\sigma-\tau)2^{-n-1},\ \ \kappa_n:=\kappa\left(1-2^{1-n}\right)
$$
and
$$
w_n:=(u-\kappa_n)^+,\ \ \Gamma_n:=(-\tau^2_n,\tau^2_n)\times B_{\tau_n}, \tilde\Gamma_n:=(-\tilde\tau^2_n,\tilde\tau^2_n)\times B_{\tilde\tau_n}.
$$
Clearly,
$$
\kappa_n\uparrow\kappa,\ \ \Gamma_{n+1}\subset \tilde\Gamma_n\subset \Gamma_n\downarrow [-\tau^2,\tau^2]\times \bar B_\tau=\bar Q_\tau.
$$
Since $\kappa_{n+1}-\kappa_n=\kappa 2^{-n}$, for any $r,s\in[1,\infty]$, by \eqref{LK1} we have 
\begin{align}\label{LK11}
\|\1_{\{w_{n+1}\not=0\}\cap\Gamma_n}\|_{\mL^{s,r}_{t,x}}\leq 2^{n}\kappa^{-1}\|\1_{\Gamma_n}w_n\|_{\mL^{s,r}_{t,x}},
\end{align}
and by \eqref{LK4}, for any $(r,s)\in\sI_\varkappa$, there is a $\delta\in(0,1)$ such that
\begin{align}\label{LK33}
\|\1_{\Gamma_{n+1}}w_{n+1}\|_{\mL^{s,r}_{t,x}}\lesssim
\tfrac{2^{n}\|w_{n+1}\|_{\sV_{\tilde\Gamma_n}}}{\sigma-\tau}\cdot\Big(2^n\kappa^{-1}\|\1_{\Gamma_n}w_n\|_{\mL^{1,1}_{t,x}}\Big)^\delta.
\end{align}
Now let $(r_i, s_i), i=1,2,3$ be as in Lemma \ref{Le32}. If we define
\begin{align*}
\ell^{(i)}_n:=\|\1_{\Gamma_n}w_{n}\|_{\mL^{s_i,r_i}_{t,x}},\ i=1,2,3,
\end{align*}
then by \eqref{LK33},  for some $\delta_i=\delta_i(r_i,s_i,\varkappa)\in(0,1)$ and $C=C(r_i,s_i,\varkappa,d)>0$,
$$
\ell^{(i)}_{n+1}\lesssim_C \tfrac{2^{n}\|w_{n+1}\|_{\sV_{\tilde\Gamma_n}}}{\sigma-\tau}\cdot\big(2^n\kappa^{-1}\|\1_{\Gamma_n}w_n\|_{\mL^{1,1}_{t,x}}\big)^{\delta_i},\ i=1,2,3.
$$
In particular, we have
\begin{align}\label{Es11}
a_{n+1}&:=\frac{1}{\kappa}\sum_{i=1}^3\ell^{(i)}_{n+1}
\lesssim \tfrac{2^{n}\|w_{n+1}\|_{\sV_{\tilde\Gamma_n}}}{(\sigma-\tau)\kappa}\sum_{i=1}^3
\big(2^n\kappa^{-1}\|\1_{\Gamma_n}w_n\|_{\mL^{1,1}_{t,x}}\big)^{\delta_i}\no
\\&\lesssim\tfrac{4^{n}\|w_{n+1}\|_{\sV_{\tilde\Gamma_n}}}{(\sigma-\tau)\kappa}\sum_{i=1}^3
\left(\frac{\ell^{(1)}_n}{\kappa}\right)^{\delta_i}\leq
 \tfrac{4^{n}\|w_{n+1}\|_{\sV_{\tilde\Gamma_n}}}{(\sigma-\tau)\kappa}\sum_{i=1}^3a_n^{\delta_i},
\end{align}
where the second inequality is due to $\|\1_{\Gamma_n}w_n\|_{\mL^{1,1}_{t,x}}\leq C\|\1_{\Gamma_n}w_n\|_{\mL^{s_1,r_1}_{t,x}}$.

On the other hand, note that
$$
0\leq w_{n+1}\leq w_n\Rightarrow |\nabla w_{n+1}|=|\nabla u|\1_{\{w_{n+1}\not=0\}}\leq|\nabla u|\1_{\{w_n\not=0\}}=|\nabla w_n|.
$$
By Lemma \ref{Le32} with $w=w_{n+1}$, $Q_1=\tilde\Gamma_n$, $Q_2=\Gamma_n$ and \eqref{LK1}, we have for some $\gamma\geq 1$,
\begin{align*}
\|w_{n+1}\|_{\sV_{\tilde\Gamma_n}}
&\lesssim_C 
2^{\gamma n}\sum_{i=1,2}\|\1_{\Gamma_n}w_{n+1}\|_{\mL^{s_i,r_i}_{t,x}}
+\|f\1_{\Gamma_n}\|_{\mL^{q_4,p_4}_{t,x}}\|\1_{\{w_{n+1}\not=0\}\cap \Gamma_n}\|_{\mL^{s_3,r_3}_{t,x}}\no\\
&\lesssim_C 2^{\gamma n}(\ell_n^{(1)}+\ell_n^{(2)})+\|f\1_{Q_2}\|_{\mL^{q_4,p_4}_{t,x}}(2^{n}\kappa^{-1}\ell_n^{(4)}),
\end{align*}
where $C=C(\Theta)>0$. This implies that for $\kappa\geq \|f\1_{Q_2}\|_{\mL^{q_4,p_4}_{t,x}}$,
$$
\|w_{n+1}\|_{\sV_{\tilde\Gamma_n}}\lesssim_C2^{\gamma n}(\ell_n^{(1)}
+\ell_n^{(2)})+2^{n}\ell_n^{(3)}\leq 2^{\gamma n}\sum_{i=1}^3\ell_n^{(i)}=2^{\gamma n}a_n\kappa .
$$
Substituting this into \eqref{Es11}, we obtain that for some $C_0,\gamma_0>1$,
$$
a_{n+1}\leq \frac{C_0 2^{\gamma_0 n}a_n}{\sigma-\tau}\sum_{i=1}^3a_n^{\delta_i},\ \ \forall n\in\mN.
$$
Let $\delta:=\delta_1\wedge\cdots\wedge\delta_4$. Suppose that
$$
\kappa\geq \left(\left[\frac{4C_02^{\gamma_0(1+\delta)/\delta}}{\sigma-\tau}\right]^{\frac1\delta}\sum_{i=1,2,3}
\|u^+\1_{Q_\sigma}\|_{\mL^{s_i,r_i}_{t,x}}\right)\vee\|f\1_{Q_2}\|_{\mL^{q_4,p_4}_{t,x}}.
$$
Then $a_1\leq (4C_02^{\gamma_0(1+\delta)/\delta})^{-\frac1\delta}$, and by Fatou's lemma and Lemma \ref{Le15},
\begin{align*}
\|(u-\kappa)^+\1_{Q_\tau}\|_{\mL^{s_1,r_1}_{t,x}}\leq \liminf_{n\to\infty}\|w_n\1_{\Gamma_n}\|_{\mL^{s_1,r_1}_{t,x}}
=\liminf_{n\to\infty}\ell^{(1)}_n\leq\kappa\cdot\limsup_{n\to\infty}a_n=0,
\end{align*}
which in turn implies that $(u-\kappa)^+=0$ on $Q_\tau$, and so,
$$
\|u^+\1_{Q_\tau}\|_\infty\leq\left(\left[\frac{4C_02^{\gamma_0(1+\delta)/\delta}}{\sigma-\tau}\right]^{\frac1\delta}
\sum_{i=1,2,3}\|u^+\1_{Q_\sigma}\|_{\mL^{s_i,r_i}_{t,x}}\right)\vee\|f\1_{Q_2}\|_{\mL^{q_4,p_4}_{t,x}}.
$$
To show \eqref{JK2}, without loss of generality, we may assume 
$$
p\leq\gamma/2,\ \ \gamma:=\max_{i=1,2,3}\{s_i, r_i\}.
$$
Thus by H\"older's inequality and Young's inequality, we have
\begin{align*}
\|u^+\1_{Q_\tau}\|_\infty&\leq C(\sigma-\tau)^{-\frac1\delta}\|u^+\1_{Q_\sigma}\|_{\mL^{\gamma,\gamma}_{t,x}}+\|f\1_{Q_2}\|_{\mL^{q_4,p_4}_{t,x}}\\
&\leq C(\sigma-\tau)^{-\frac1\delta}\|u^+\1_{Q_2}\|^{1-\frac{p}{\gamma}}_\infty
\|u^+\1_{Q_\sigma}\|^{\frac{p}{\gamma}}_{\mL^{p,p}_{t,x}}+\|f\1_{Q_2}\|_{\mL^{q_4,p_4}_{t,x}}\\
&\leq\tfrac12\|u^+\1_{Q_\sigma}\|_\infty+C(\sigma-\tau)^{-\frac{\gamma}{p\delta}}\|u^+\1_{Q_2}\|_{\mL^{p,p}_{t,x}}+\|f\1_{Q_2}\|_{\mL^{q_4,p_4}_{t,x}},
\end{align*}
where $C=C(\Theta)$ is independent of $\sigma,\tau$.
By Lemma \ref{Le26}, we conclude the proof.
\end{proof}

\subsection{Proof of Theorem \ref{TH22}}
Without loss of generality, we assume $T=1$ and 
$$
u(t,x)=f(t,x)\equiv 0, \ \  \forall t\leq 0. 
$$
For $z\in\mR^d$, we write
$$
Q^z_i:=Q^{0,z}_i,\ \ i=1,2,3.
$$
Let $u$ be a Lipschitz weak solution of PDE \eqref{PDE0} in the sense of Definition \ref{Def1}.
By translation and Lemma \ref{Le32} with $w=u^+, u^-$,  there is a constant $C=C(\Theta)>0$ such that for all $t\in[0,1]$,
$$
\|u\cI^t\|_{\sV_{Q^z_1}}\lesssim_C \sum_{i=1,2}\|\1_{Q^z_2}u\cI^t\|_{\mL^{s_i,r_i}_{t,x}}
+\|f\1_{Q^z_2}\cI^t\|_{\mL^{q_4,p_4}_{t,x}},
$$
where $(r_i,s_i)$ are as in Lemma \ref{Le32}.
By Lemma \ref{Le22}, we further have for some $\beta\in(1,\infty)$ and any $\eps\in(0,1)$,
\begin{align*}
\|u\cI^t\|_{\sV_{Q^z_1}}&\leq\eps\|\1_{Q^z_3}\nabla u\cI^t\|_{\mL^{\varkappa,2}_{x,t}}+C_\eps\|\1_{Q^z_3}u\cI^t\|_{\mL^{\beta,2}_{t,x}}
+C\|f\1_{Q^z_2}\cI^t\|_{\mL^{q_4,p_4}_{t,x}}.
\end{align*}
Taking supremum in $z\in\mR^d$ for both sides and by \eqref{DK1}, we obtain
$$
\nor u\cI^t\nor_{\widetilde\sV}\lesssim\sup_{z}\|u\cI^t\|_{\sV_{Q^z_1}}\lesssim\eps
\nor\nabla u\cI^t\nor_{\widetilde\mL^{\varkappa,2}_{x,t}}+\nor u\cI^t\nor_{\widetilde\mL^{\beta,2}_{t,x}}
+\nor  f\cI^t\nor_{\mL^{q_4,p_4}_{t,x}}.
$$
Since $\nor\nabla u\cI^t\nor_{\widetilde\mL^{\varkappa,2}_{x,t}}\leq \nor u\cI^t\nor_{\widetilde\sV}$, choosing $\eps$ small enough, we obtain
\begin{align}\label{Es151}
\nor u\cI^t\nor_{\widetilde\sV}\lesssim_C\nor u\cI^t\nor_{\widetilde\mL^{\beta,2}_{t,x}}+\nor  f\cI^t\nor_{\mL^{q_4,p_4}_{t,x}}.
\end{align}
Since $\beta<\infty$ and $u(t)\equiv 0$ for $t\leq 0$,  the above inequality implies that for any $t\in[0,1]$,
$$
\nor u(t)\nor_{2}\lesssim_C\left(\int^t_0\nor u(s)\nor^{\beta}_{2}\dif s\right)^{1/\beta}+\nor  f\cI^t\nor_{\mL^{q_4,p_4}_{t,x}}.
$$
By Gronwall's inequality we obtain 
$$
\nor u\cI^t\nor_{\widetilde\mL^{\infty,2}_{t,x}}\leq\sup_{t\in[0,1]}\nor u(t)\nor_{2}\lesssim \nor  f\cI^t\nor_{\mL^{q_4,p_4}_{t,x}}.
$$
which together with \eqref{Es151} yields 
\begin{align}\label{DP2}
\nor u\1_{[0,1]}\nor_{\widetilde\sV}\lesssim_C\nor  f\1_{[0,1]}\nor_{\mL^{q_4,p_4}_{t,x}}.
\end{align}
Finally,  by \eqref{JK2} and \eqref{DP2}, we also have
\begin{align*}
&\|u\|_{L^\infty([0,1]\times\mR^d)}\leq \sup_{z}\|(u^++u^{-})\1_{[0,1]\times B^z_1}\|_{\infty}\\
&\quad\lesssim \nor u\1_{[0,1]}\nor_{\widetilde\mL^{2,2}_{t,x}}+\nor  f\1_{[0,1]}\nor_{\mL^{q_4,p_4}_{t,x}} \lesssim \nor  f\1_{[0,1]}\nor_{\mL^{q_4,p_4}_{t,x}}.
\end{align*}
The proof is complete.

\section{Weak solutions of SDEs with rough coefficients}\label{Se4}

In this section we present an application of the global boundedness \eqref{Es161} in SDEs, and show the existence of weak solutions to 
SDE \eqref{SDE0} with rough coefficients.
First of all, we recall the following notion of weak solutions to SDE \eqref{SDE0}.
\bd\label{Def41}
Let $\frak{F}:=(\Omega,\sF,\bP; (\sF_t)_{t\geq 0})$ be a stochastic basis  and $(X,W)$ a pair of $\sF_t$-adapted processes defined thereon. 
Given $(s,x)\in\mR_+\times\mR^d$, we call triple $(\frak{F}, X,W)$ a weak solution of SDE \eqref{SDE0} with starting point $x\in\mR^d$
at time $s$ if
\begin{enumerate}[(i)]
\item $\bP(X_t=x, t\in[0,s])=1$ and $W$ is an $\sF_t$-Brownian motion;
\item for all $t\geq s$, it holds that $\bP$-a.s.,
$$
\int^t_s\big(|\sigma(r,X_r)|^2+|b(r,X_r)|\big)\dif r<\infty,
$$
and
\begin{align*}
X_t=x+\sqrt{2}\int^t_s\sigma(r,X_r)\dif W_r+\int^t_sb(r,X_r)\dif r.
\end{align*}
\end{enumerate}
\ed

Recall $p_0, p_1$ from \eqref{PP0} and the convention that the repeated indices will be summed automatically, for instances,
$$
\p_ia^{ij}=\sum_{i=1}^{d}\p_ia^{ij}, \quad \p_i\p_ja^{ij}=\sum_{i,j=1}^{d}\p_i\p_ja^{ij}.
$$
We introduce the following assumptions on $\sigma$ and $b$:
\begin{enumerate}[{\bf ($\widetilde {\bf H}^\sigma$)}]
\item  Suppose that there are a sequence of $d\times d$-matrix functions $\sigma_n\in L^\infty(\mR_+; C^\infty_b)$, $(p_2,q_2)\in\mI^d_{p_0}$
  and $\kappa_0>0$ such that for all $n\in\mN$,
\begin{align}\label{NA0}
\,\,\qquad\nor\lambda^{-1}_n\nor_{p_0}+\nor\mu_n\nor_{p_1}+\nor \p_ia^{ij}_n\nor_{\widetilde\mL^{p_1,\infty}_{x,t}}
+\nor (\p_i\p_ja^{ij}_n)^+\nor_{\widetilde\mL^{q_2,p_2}_{t,x}}\leq\kappa_0,
\end{align}
where $a_n:=\sigma_n\sigma^*_n$, $\lambda_n$ and $\mu_n$ are defined as in \eqref{DW3} by $a_n$. Moreover, 
for some $p_4,q_4\in[2,\infty]$ with $(\frac{p_4}{2},\frac{q_4}{2})\in\mI^d_{p_0}$ and for any $T,R>0$,
\begin{align}\label{NA1}
\qquad\sup_n\nor\sigma_n\nor_{\widetilde\mL^{q_4,p_4}_{t,x}}=:\kappa_1<\infty,\quad 
\lim_{n\to\infty}\|(\sigma_n-\sigma)\1_{[0,T]\times B_R}\|_{\mL^{q_4,p_4}_{t,x}}=0.
\end{align}
\end{enumerate}

\begin{enumerate}[{\bf ($\widetilde {\bf H}^b$)}]
\item Let $b=b_1+b_2$ satisfy {\bf (${\bf H}^b$)} and belong to $\widetilde\mL^{q_4,p_4}_{t,x}$  for some $(p_4,q_4)\in\mI^d_{p_0}$. 
\end{enumerate}

Let $\Theta$ be defined by \eqref{TH}. Below we shall write
$$
\widetilde\Theta:=\Big(\Theta, p_i,q_i,\kappa_0,\kappa_1, \nor b\nor_{\widetilde\mL^{q_4,p_4}_{t,x}}\Big).
$$
We have the following existence result.
\bt\label{Th15}
Under {\bf ($\widetilde {\bf H}^\sigma$)} and {\bf ($\widetilde {\bf H}^b$)}, for any $(s,x)\in\mR_+\times\mR^d$, there is a weak solution 
$(\frak{F}, X,W)$ for SDE \eqref{SDE0} starting from $x$ at time $s$. 
Moreover, for any $(p,q)\in\mI^d_{p_0}$ and $T>s$, there are $\theta\in(0,1)$ and constant 
$C=C(T,\widetilde\Theta,p,q)>0$ such that for any $s\leq t_0<t_1\leq T$ and $f\in\widetilde\mL^{q,p}_{t,x}$, 
\begin{align}\label{DD87}
\bE\left(\int^{t_1}_{t_0}f(r, X_r)\dif r\Big|\sF_{t_0}\right)\leq C(t_1-t_0)^\theta\nor f\nor_{\widetilde\mL^{q,p}_{t,x}}.
\end{align}
\et

In the following proof, we assume $s=0$ and $x\in\mR^d$.
Let $\sigma_n$ be as in {\bf ($\widetilde {\bf H}^\sigma$)} and $b_n(t,x)=b*\rho_n(t,x)$ be the mollifying approximation of $b$. 
In particular, 
\begin{align}\label{CC1}
\sigma_n, b_n\in L^\infty([0,T];C^\infty_b(\mR^d)),
\end{align}
and the following SDE admits a unique strong solution (cf. \cite{St-Va}):
\begin{align}\label{SDE9}
\dif X^n_t=b_n(t,X^n_t)\dif t+\sqrt{2}\sigma_n(t,X^n_t)\dif W_t,\ \ X^n_0=x.
\end{align}
Note that the generator of SDE \eqref{SDE9} is given by
\begin{align*}
\sL^{\sigma_n,b_n}_t f(x)&=(\sigma^{ik}_n\sigma^{jk}_n)(t,x)\p_i\p_j f(x)+b^j_n(t,x)\p_j f(x)\\
&=\p_i(a^{ij}_n(t,\cdot)\p_j f)(x)+\tilde b^j_n(t,x)\p_j f(x),
\end{align*}
where
$$
a^{ij}_n:=\sigma^{ik}_n\sigma^{jk}_n,\ \ \tilde b^j_n:=b^j_n-\p_ia^{ij}_n.
$$
In particular, by {\bf ($\widetilde {\bf H}^\sigma$)}, one sees that
{\bf (${\bf H}^{a}$)} holds for $a_n$ uniformly in $n$, and {\bf (${\bf H}^b$)} holds for $\tilde b_n=b_{1,n}+(b_{2,n}-\p_ia^{i\cdot}_n)$ uniformly in $n$, where $b_{i,n}:=b_i*\rho_n$.

We first show the following key Krylov estimate (see \cite{Zh-Zh}).
\bt\label{Th43}
Under {\bf ($\widetilde {\bf H}^\sigma$)} and {\bf ($\widetilde {\bf H}^b$)}, 
for any $(p,q)\in\mI^d_{p_0}$ and $T>0$, there are $\theta=\theta(p,q,d,p_0)\in(0,1)$ and $C=C(T,\widetilde\Theta,p,q)>0$ independent of 
starting point $x$ such that
for any $0\leq t_0<t_1\leq T$ and $f\in\widetilde\mL^{q,p}_{t,x}$, 
\begin{align}\label{DD7}
\sup_n\mE\left(\int^{t_1}_{t_0}f(s, X^n_s)\dif s\Big|\sF_{t_0}\right)\leq C(t_1-t_0)^\theta\nor f\nor_{\widetilde\mL^{q,p}_{t,x}}.
\end{align}
\et
\begin{proof}
Below we drop the super and subscripts $n$ for simplicity.
Without loss of generality, we may assume $f\in C^\infty_0(\mR^{d+1})$. Fix $t_1\in(0, T]$ and consider the following backward PDEs:
$$
\p_t u+\sL^{\sigma,b}_t u=f,\ \ u(t_1)=0.
$$
Under \eqref{CC1}, it is well known that there is a unique solution $u\in L^\infty_{loc}([0,t_1]; C^\infty_b(\mR^d))$ so that (cf. \cite{St-Va})
$$
u(t,x)=\int^{t_1}_t(\sL^{\sigma,b}_s u-f)(s, x)\dif s,\ \ t\in[0,t_1].
$$
By It\^o's formula, for any $t_0\leq t_1$, we have
$$
u(t_1, X_{t_1})-u(t_0, X_{t_0})=\int^{t_1}_{t_0}f(s, X_s)\dif s+\sqrt{2}\int^{t_1}_{t_0}(\sigma^*\nabla u)(s, X_s)\dif W_s.
$$
Taking conditional expectations with respect to $\sF_{t_0}$, we obtain
\begin{align}\label{DD8}
\mE\left(\int^{t_1}_{t_0}f(s, X_s)\dif s\Big|\sF_{t_0}\right)\leq \|u\|_{L^\infty([t_0,t_1]\times\mR^d)}.
\end{align}
On the other hand, since $(p,q)\in\mI^d_{p_0}$, we can choose $q'<q$ so that
$$
\tfrac{1}{p}<\big(1-\tfrac{1}{q'}\big)\big(\tfrac{2}{d}-\tfrac{1}{{p_0}}\big).
$$
Thus by the assumptions, \eqref{Es161} of Theorem \ref{TH22}, there exists a 
constant $C=C(T,\widetilde\theta,p,q)>0$ independent of $n$ such that
$$
\|u\|_{L^\infty([t_0,t_1]\times\mR^d)}\lesssim_C\nor  f\1_{[t_0,t_1]}\nor_{\mL^{q',p}_{t,x}}
\lesssim_C (t_1-t_0)^\theta\nor f \nor_{\mL^{q,p}_{t,x}},
$$
where $\theta=\frac{1}{q'}-\frac{1}{q}$ and the second inequality is due to H\"older's inequality.
Substituting it into \eqref{DD8}, we obtain \eqref{DD7}.
\end{proof}
We need the following simple lemma.
\bl\label{Le44}
Let $(X_t)_{t\geq 0}$ be a right continuous stochastic process on a filtered probability space $(\Omega,\sF,\mP; (\sF_t)_{t\geq 0})$.
Suppose that for some $Y\in L^1(\Omega)$ and $A>0$,
$$
|X_t|\leq Y,\ \ \mE(X_t|\sF_t)\leq A,\ \ \mP-a.s.
$$
Then for any finite stopping time $\tau$, it holds that
$$
\mE(X_\tau|\sF_\tau)\leq A,\ \ \mP-a.s.
$$
\el
\begin{proof}
Let $\tau_n$ be a sequence of decreasing stopping times with values in $\mD:=\{k\cdot 2^{-n}: k, n\in \mN\}$ 
and so that $\tau_n\to\tau$ as $n\to\infty$.
Note that for each $n\in\mN$,
\begin{align*}
\mE\left(X_{\tau_n}|\sF_{\tau_n} \right)
=\mE\left(\sum_{t\in\mD}\1_{\{\tau_n=t\}}X_t|\sF_{\tau_n} \right)
=\sum_{t\in\mD}\1_{\{\tau_n=t\}}\mE\left(X_t|\sF_t \right)\leq A.
\end{align*}
By the dominated convergence theorem and $\sF_\tau\subset\sF_{\tau_n}$, we have
\begin{align*}
\mE\left(X_{\tau}|\sF_{\tau} \right)
&=\lim_{n\to\infty}\mE\left(X_{\tau_n}|\sF_{\tau} \right)
=\lim_{n\to\infty}\mE\left(X_{\tau_n}|\sF_{\tau_n}|\sF_{\tau} \right)\leq A.
\end{align*}
The proof is complete.
\end{proof}
\br\rm
By this lemma, one sees that \eqref{DD7} is equivalent to that for any stopping time $\tau\leq T$, $\delta\in(0,1)$ and $f\in\widetilde\mL^{q,p}_{t,x}$, 
\begin{align}\label{DD77}
\sup_n\mE\left(\int^{\tau+\delta}_{\tau}f(s, X^n_s)\dif s\Big|\sF_{\tau}\right)\leq C\delta^\theta\nor f\nor_{\widetilde\mL^{q,p}_{t,x}}.
\end{align}
\er
\bl\label{Le33}
Under {\bf ($\widetilde {\bf H}^\sigma$)} and {\bf ($\widetilde {\bf H}^b$)}, for any $T>0$, there are $\theta\in(0,1)$ and 
constant $C=C(T,\widetilde\Theta)>0$ such that for all $\delta\in(0,1)$,
\begin{align}\label{NB1}
\sup_n\mE\left(\sup_{t\in[0,T]}\sup_{s\in[0,\delta]}|X^n_{t+s}-X^n_{t}|^{1/2}\right)\leq C\delta^{\theta/2}
\end{align}
and
\begin{align}\label{NB2}
\sup_n\mE\left(\sup_{t\in[0,T]}|X^n_{t}|\right)\leq C.
\end{align}
\el
\begin{proof}
We only prove \eqref{NB1}. Let $\tau$ be any stopping time less than $T$. Notice that
\begin{align*}
X^n_{\tau+t}-X^n_{\tau}&=\int^{\tau+t}_\tau b_n(s,X^n_s)\dif s+\sqrt{2}\int^{\tau+t}_\tau\sigma_n(s, X^n_s)\dif W_s,\ \ t>0.
\end{align*}
By Burkholder's inequality and \eqref{DD77}, we have
\begin{align*}
\mE\left(\sup_{t\in[0,\delta]}|X^n_{\tau+t}-X^n_{\tau}|\right)&
\lesssim \mE \int^{\tau+\delta}_\tau |b_n|(s,X^n_s)\dif s+\left(\mE\int^{\tau+\delta}_\tau|\sigma_n(s, X^n_s)|^2\dif s\right)^{1/2}\\
&\leq  C\delta^\theta\nor b_n\nor_{\widetilde\mL^{q_4,p_4}_{t,x}}+C\delta^\theta\nor\sigma_n\nor_{\widetilde\mL^{q_4,p_4}_{t,x}}
\leq C\delta^\theta,
\end{align*}
where $C$ is independent of $n$ and $\delta$.
Thus by \cite[Lemma 2.7]{Zh-Zh1}, we obtain \eqref{NB1}.
\end{proof}

Let $\mC$ be the space of all continuous functions from $\mR_+$ to $\mR^d$, which is endowed with the locally uniformly metric 
so that $\mC$ becomes a Polish space. 
Let $\mQ_n$ be the law of $(X^n_\cdot, W_\cdot)$ in product space $\mC\times\mC$. 
For each $x\in\mR^d$, by Lemma \ref{Le33} and \cite[Theorem 1.3.2]{St-Va}, the law of $X^n_\cdot$ is tight in $\mC$.
By Lemma \ref{Le33} and Prokhorov's theorem, there are a subsequence still denoted by $n$ and $\mQ\in\cP(\mC\times\mC)$ so that
$$
\mQ_n\to \mQ\mbox{ weakly.}
$$
Now, by Skorokhod's representation theorem,  there are a probability space $(\tilde\Omega,\tilde\sF,\tilde\mP)$ and 
random variables $(\tilde X^n , \tilde W^n)$ and $(\tilde X,\tilde W)$ defined on it such that
\begin{align}\label{FD4}
(\tilde X^n , \tilde W^n)\to (\tilde X,\tilde W),\ \ \tilde\mP-a.s.
\end{align}
and
\begin{align}\label{FD5}
\tilde\mP\circ(\tilde X^n , \tilde W^n)^{-1}=\mQ_n=\mP\circ(X^n , W)^{-1},\quad
\tilde\mP\circ(\tilde X, \tilde W)^{-1}=\mQ.
\end{align}
Define $\tilde\sF^n _t:=\sigma(\tilde W^n _s, \tilde X^n_s;s\leq t)$. 
Notice that
$$
\mP(W _t-W_s\in\cdot |\sF _s)=\mP(W _t-W _s\in\cdot)
$$
implies that
$$
\tilde\mP(\tilde W^n _t-\tilde W^n _s\in\cdot |\tilde \sF^n _s)=\tilde\mP(\tilde W^n _t-\tilde W^n _s\in\cdot).
$$
In other words, $\tilde W^n_t$ is an $\tilde\sF_t^n $-Brownian motion. Thus, by \eqref{SDE9} and \eqref{FD5} we have
\begin{align}\label{KK1}
\tilde X^n _t=x+\int^t_0b_n\big(s,\tilde X^n_s\big)\dif s+\int^t_0\sigma_n\big(s,\tilde X^n_s\big)\dif \tilde W^n_s.
\end{align}
Moreover, by \eqref{DD7}, we also have
\begin{align}\label{DD97}
\sup_n\tilde\mE\left(\int^{t_1}_{t_0}f(s, \tilde X^n_s)\dif s\Big|\tilde\sF_{t_0}\right)\leq C(t_1-t_0)^\theta\nor f\nor_{\widetilde\mL^{q,p}_{t,x}}.
\end{align}

In order to take the limits, we recall a result of Skorokhod \cite[p.32]{Sk}. 
\bl\label{Le332}
Let $\{f_n(t),t\geq 0, n\in\mN\}$ be a sequence of measurable $\tilde\sF^n_t$-adapted processes. Suppose that
for every $T, \eps>0$, there is an $M_\eps>0$ such that
$$
\sup_n\tilde\mP\left\{\sup_{t\in[0,T]}|f_n(t)|>M_\eps\right\}\leq\eps,
$$
and also,
$$
\lim_{n\to\infty}\tilde\mP\left\{\sup_{t\in[0,T]}|f_n(t)-f(t)|>\eps\right\}=0.
$$
Then it holds that for every $T>0$, 
$$
\int^T_0 f_n(t)\dif \tilde W^n_t\to\int^T_0 f(t)\dif \tilde W_t\mbox{ in probability as $n\to\infty$}.
$$
\el

\bl\label{Le302}
For each $t>0$, the following limits hold in probability as $n\to \infty$,
\begin{align}
\int^t_0b_n\big(s, \tilde X^n _s\big)\dif s
&\to \int^t_0b\big(s,\tilde X_{s}\big)\dif s,\label{Lim1}\\
\int^t_0\sigma_n\big(s, \tilde X^n _s\big)\dif \tilde W^n _s
&\to \int^t_0\sigma\big(s,\tilde X_{s}\big)\dif \tilde W_s.\label{Lim2}
\end{align}
\el
\begin{proof}
We only prove \eqref{Lim2}. For simplicity, we shall write
$\sigma_\infty:=\sigma$ and drop the tilde.  For each $n\in\mN_\infty=\mN\cup\{\infty\}$,
let $\sigma^\eps_n(t,x):=\sigma_n*\rho_\eps(t,x)$ be the mollifying approximation of $\sigma_n$. 
It suffices to show the following two limits:  for fixed $\eps>0$,
\begin{align}\label{Lim4}
\int^t_0\sigma^\eps_n\big(s,  X^n _s\big)\dif  W^n _s
\to\int^t_0\sigma^\eps_\infty\big(s, X_{s}\big)\dif  W_s\mbox{ in probability $n\to\infty$},
\end{align}
and
\begin{align}\label{Lim3}
\lim_{\eps\to 0}\sup_{n\in\mN_\infty}\mE\left|\int^t_0(\sigma^\eps_n\big(s,  X^n _s\big)-\sigma_n\big(s,  X^n _s\big))\dif  W^n _s\right|^2=0.
\end{align}
Clearly, limit \eqref{Lim4} follows by Lemma \ref{Le332}. We look at \eqref{Lim3}.
For $R>0$, we define
$$
\tau^n_R:=\inf\{t>0: |X^n_t|\geq R\}.
$$
By \eqref{NB2}, we have
\begin{align}\label{FS7}
\lim_{R\to\infty}\sup_n\mP(\tau^n_R\leq t)\leq \lim_{R\to\infty}\sup_n\frac{1}{R}\mE\left(\sup_{s\in[0,t]}|X^n_s|\right)=0.
\end{align}
For \eqref{Lim3}, by It\^o's isometry we have
\begin{align}\label{JL1}
\begin{split}
&\mE\left|\int^t_0(\sigma^\eps_n\big(s,  X^n _s\big)-\sigma_n\big(s,  X^n _s\big))\dif  W^n _s\right|^2
\\&\qquad\leq\mE\int^t_0|\sigma^\eps_n\big(s,  X^n _s\big)-\sigma_n\big(s,  X^n _s\big)|^2\dif s=:I^n_R(\eps)+J^n_R(\eps),
\end{split}
\end{align}
where
\begin{align*}
I^n_R(\eps)&:=\mE\left(\1_{\{\tau^n_R>t\}}\int^t_0|\sigma^\eps_n\big(s,  X^n _s\big)-\sigma_n\big(s,  X^n _s\big)|^2\dif s\right),\\
J^n_R(\eps)&:=\mE\left(\1_{\{\tau^n_R\leq t\}}\int^{t}_0|\sigma^\eps_n\big(s,  X^n _s\big)-\sigma_n\big(s,  X^n _s\big)|^2\dif s\right).
\end{align*}
For $I^n_R(\eps)$, since $(\frac{p_4}2,\frac{q_4}2)\in\mI^d_{p_0}$, by \eqref{DD97} we have
\begin{align*}
I^n_R(\eps)&\leq\mE\left(\int^t_0\1_{| X^n_s|\leq R}|\sigma^\eps_n\big(s,  X^n _s\big)-\sigma_n\big(s,  X^n _s\big)|^2\dif s\right)\\
&\lesssim\|\1_{[0,t]\times B_R}(\sigma^\eps_n-\sigma_n)\|^2_{\mL^{q_4,p_4}_{t,x}},
\end{align*}
where the implicit constant is independent of $n,\eps,R$.
For each $R>0$, since 
$$
\lim_{n\to\infty}\sup_{\eps\in(0,1)}\|\1_{[0,t]\times B_R}(\sigma^\eps_n-\sigma^\eps)\|_{\mL^{q_4,p_4}_{t,x}}
\leq\lim_{n\to\infty}\|\1_{[0,t]\times B_{2R}}(\sigma_n-\sigma)\|_{\mL^{q_4,p_4}_{t,x}}=0,
$$ 
and for each $n\in\mN_\infty$,
$$
\lim_{\eps\to 0}\|\1_{[0,t]\times B_R}(\sigma^\eps_n-\sigma_n)\|_{\mL^{q_4,p_4}_{t,x}}=0,
$$
it follows that for each $R>0$,
\begin{align}\label{JL2}
\lim_{\eps\to 0}\sup_n I^n_R(\eps)\lesssim\lim_{\eps\to 0}\sup_n\|\1_{[0,t]\times B_R}(\sigma^\eps_n-\sigma_n)\|^2_{\mL^{q_4,p_4}_{t,x}}=0.
\end{align}
For $J^n_R(\eps)$, since $(\frac{p_4}{2},\frac{q_4}{2})\in\mI^d_{p_0}$, one can choose $\gamma>1$ being close to $1$ 
so that $(\frac{p_4}{2\gamma},\frac{q_4}{2\gamma})\in\mI^d_{p_0}$.
By H\"older's inequality and \eqref{DD7} we have
\begin{align*}
J^n_R(\eps)&\leq(\mP(\tau^n_R\leq t))^{\frac{\gamma-1}{\gamma}}\left(\mE\int^{t}_0
|\sigma^\eps_n\big(s,  X^n _s\big)-\sigma_n\big(s,  X^n _s\big)|^{2\gamma}\dif s\right)^{\frac1\gamma}\\
&\lesssim (\mP(\tau^n_R\leq t))^{\frac{\gamma-1}{\gamma}}\nor\sigma^\eps_n-\sigma_n\nor^{2\gamma}_{\mL^{q_4,p_4}_{t,x}}
\lesssim (\mP(\tau^n_R\leq t))^{\frac{\gamma-1}{\gamma}},
\end{align*}
where the implicit constant is independent of $\eps,n,R$.
By \eqref{FS7}, we have
\begin{align}\label{JL3}
\lim_{R\to\infty}\sup_n\sup_{\eps} J^n_R(\eps)=0.
\end{align}
Combining \eqref{JL1}, \eqref{JL2} and \eqref{JL3}, we obtain \eqref{Lim3}. The proof is complete.
\end{proof}
\begin{proof}[Proof of Theorem \ref{Th15}]
It follows by taking limits for both sides of \eqref{KK1} and Lemma \ref{Le302}.
As for \eqref{DD87}, it follows by taking limits for \eqref{DD97} with $f\in C_0(\mR_+\times\mR^d)$.
\end{proof}

\section{Strong Markov selection}
In this section we use Krylov's Markov selection theorem to show the existence of a strong Markov solution under 
{\bf ($\widetilde {\bf H}^\sigma$)} and {\bf ($\widetilde {\bf H}^b$)}. 
Let $\omega_t$ be the coordinate process on the continuous function space $\mC$ and 
$\cB_t:=\sigma\{\omega_s: s\leq t\}$ the natural $\sigma$-filtration.
We first recall the following notion of local martingale solutions to SDE \eqref{SDE0}.
\bd
Let $(s,x)\in\mR_+\times\mR^d$. A probability measure $\mP\in\cP(\mC)$ is called a local martingale solution of SDE \eqref{SDE0} 
starting from $x$ at time $s$ if
\begin{enumerate}[(i)]
\item $\mP(\omega_t=x, t\in[0,s])=1$ and for each $t>s$,
$$
\mP\left(\omega: \int^t_s(|b(r,\omega_r)|+|(\sigma\sigma^*)(r,\omega_r)|)\dif s<\infty\right)=1.
$$

\item For any $f\in C^\infty(\mR^d)$, the process
$$
M^f_t(\omega):=f(\omega_t)-f(\omega_s)-\int^t_s\sL^{\sigma,b}_r f(\omega_r)\dif r
$$
is a continuous local $\cB_t$-martingale after time $s$.
\end{enumerate}
The set of all the local martingale solutions of \eqref{SDE0} is denoted by $\cM^{\sigma,b}_{s,x}\subset\cP(\mC)$.
\ed

By It\^o's formula, it is easy to see that the law of a weak solution in Definition \ref{Def41} is a local martingale solution.
Moreover, we also have the following opposite conclusion (see \cite[p314, Proposition 4.11]{Ka-Sh}).
\bt
For any $\mP\in\cM^{\sigma,b}_{s,x}$, there is a weak solution $(\frak{F}, X,W)$ starting from $x$ at time $s$,
where $\frak{F}=(\Omega,\sF,\bP; (\sF_t)_{t\geq 0})$ is a stochastic basis, and so that
$$
\mP=\bP\circ X^{-1}.
$$
\et
One also needs the following notion about Krylov's estimate (see Theorem \ref{Th43}).
\bd
Let $p,q\in[1,\infty)$ and $s\geq 0$. We call a probability measure $\mP\in\cP(\mC)$ satisfy the Krylov estimate with indices $p,q$ and starting from $s$ if for any $T>s$,
there are constants $\kappa,\theta>0$ such that  for any $s\leq t_0<t_1\leq T$ and $f\in C_0(\mR_+\times\mR^d)$, 
\begin{align}\label{DD90}
\mE^\mP\left(\int^{t_1}_{t_0} f(r,\omega_r)\dif r\Big|\cB_{t_0}\right)\leq \kappa(t_1-t_0)^\theta\| f\|_{\mL^{q,p}_{t,x}}.
\end{align}
We shall denote by $\cK^{p,q}_{s,T}$ the set of all the above $\mP$.
\ed
\br\rm
By a standard approximation, \eqref{DD90} holds for all $f\in\mL^{q,p}_{t,x}$.
\er

Now we show the following main result of this section.
\bt\label{Th55}
Assume {\bf ($\widetilde {\bf H}^\sigma$)} and {\bf ($\widetilde {\bf H}^b$)}. For given $(s,x)\in\mR_+\times\mR^d$, let 
$$
\sC(s,x):=\cap_{(p,q)\in\mI^d_{p_0}, T>s}\cK^{p,q}_{s, T}\cap \cM^{\sigma,b}_{s,x}.
$$
Then $\sC(s,x)$ is a non-empty convex subset of $\cP(\mC)$ and satisfies {\bf (C1)}, {\bf (C2)} and {\bf (C3)} in appendix. 
In particular, there is a measurable mapping
$$
\mR_+\times\mR^d\ni (s,x)\mapsto \mP_{s,x}\in\sC(s,x)
$$ 
so that for each fixed $(s,x)\in\mR_+\times\mR^d$ and finite stopping time $\tau\geq s$, 
there is a $\mP_{s,x}$-null set $N\in\cB_\tau$ such that for all $\omega\notin N$,
$$
\mP_{s,x}(\cdot|\cB_\tau)(\omega)=\mP_{\tau(\omega), \omega_{\tau(\omega)}}(\cdot).
$$
\et
\begin{proof}
First of all, by Theorem \ref{Th43}, for each $(s,x)$, $\sC(s,x)$ is non-empty and convex, and for given $(p,q)\in\mI^d_{p_0}$, 
the constants $\kappa,\theta$ appearing in \eqref{DD90} are independent of $s,x$. 

{\bf Verification of (C1)} Let $(s_n,x_n)$ converge to $(s,x)$ and $\mP^n\in\sC(s_n,x_n)$. We want to show that $(\mP^n)_{n\in\mN}$ is tight.
By the equivalence between martingale solutions and weak solutions, for each $n\in\mN$, there exists a weak solution  $(\frak{F}^n, X^n,W^n)$
starting from $x_n$ at time $s_n$, where
$\frak{F}_n:=(\Omega^n,\sF^n,\bP^n; (\sF^n_t)_{t\geq 0})$, so that
$$
\mP^n=\bP^n\circ (X^n)^{-1}.
$$
Note that
$$
X^n_t=x_n+\sqrt{2}\int^t_{s_n}\sigma(r,X^n_r)\dif W^n_r+\int^t_{s_n}b(r,X^n_r)\dif r,\ \ t>s_n.
$$
Since $\mP^n\in \cap_{(p,q)\in\mI_{p_0}, T>s_n}\cK^{p,q}_{s_n,T}$, 
and the constants $\kappa,\theta$ appearing in \eqref{DD90} are independent of $n$, 
as in Lemma \ref{Le33}, one can show that for each $T>\sup s_n+1$,
$$
\sup_n\bE_n\left(\sup_{t\in[0,T]}\sup_{s\in[0,\delta]}|X^n_{t+s}-X^n_{t}|^{1/2}\right)\leq C\delta^{\theta/2},\ \ \delta\in(0,1),
$$
where $\bE_n$ stands for the expectation with respect to $\bP^n$.
Thus $(\mP^n)_{n\in\mN}$ is tight.
Let $\mP$ be any accumulation point of $\mP^n$. If necessary, by substracting a subsequence, without loss of generality we assume $\mP^n$
weakly converges to $\mP$. 
For given compact support continuous function $f$, by taking weak limits for
$$
\mE^{\mP_n}\left(\int^{t_1}_{t_0} f(r,\omega_r)\dif r\Big|\cB_{t_0}\right)\leq \kappa(t_1-t_0)^\theta\| f\|_{\mL^{q,p}_{t,x}},
$$
one sees that
$$
\mP\in \cap_{(p,q)\in\mI_{p_0}, T>s}\cK^{p,q}_{s,T}.
$$
Moreover, as in the proof in Section \ref{Se4}, one can show that $\mP\in\cM^{\sigma,b}_{s,x}$.

{\bf Verification of (C2)} Let $\mP\in\sC(s,x)$ and $\tau\geq s$ be a finite stopping time. 
Let $Q_\omega$ be a r.c.p.d. of $\mP(\cdot|\cB_\tau)$. 
By \cite[Theorem 6.1.3]{St-Va}, there is a $\mP$-null set $N_1\in\cB_\tau$ such that for all $\omega\notin N_1$,
\begin{align}\label{DA0}
Q_\omega\in\cM^{\sigma,b}_{\tau(\omega),\omega_{\tau(\omega)}}.
\end{align}
On the other hand, for fixed $p,q\in\mI^d_{p_0}$, $\delta\in(0,1)$ and $T>s+\delta$, since $\mP\in\cK^{p,q}_{s, T}$, we have
for all $t\in[s,T-\delta]$ and $f\in C_0(\mR_+\times\mR^d)$,
$$
\mE^{\mP}\left(\int^{t+\delta}_t f(r,\omega_r)\dif r\Big|\cB_{t}\right)
\leq \kappa\delta^\theta\| f\|_{\mL^{q,p}_{t,x}},\ \mP-a.s.
$$
By Lemma \ref{Le74}, there is a $\mP$-null set $N=N(p,q,f,T)\in\cB_\tau$ such that for all $\omega\notin N$ and $T>\tau(\omega)+\delta$,
$t\in[\tau(\omega),T-\delta]$,
$$
\mE^{Q_\omega}\left(\int^{t+\delta}_t f(r,\omega_r)\dif r\Big|\cB_{t}\right)
\leq \kappa\delta^\theta\| f\|_{\mL^{q,p}_{t,x}},\ Q_\omega-a.s.,
$$
Since $C_0(\mR_+\times\mR^d)$ and $\mI^d_{p_0}$ are separable, one can find a common $\mP$-null set $N_2$
such that for all $\omega\notin N_2$ and $(p,q)\in\mI^d_{p_0}$, $\delta\in(0,1)$, $T>\tau(\omega)+\delta$, $t\in[\tau(\omega),T-\delta)$,
$f\in C_0(\mR_+\times\mR^d)$, 
$$
\mE^{Q_\omega}\left(\int^{t+\delta}_t f(r,\omega_r)\dif r\Big|\cB_{t}\right)
\leq \kappa\delta^\theta\| f\|_{\mL^{q,p}_{t,x}},\ Q_\omega-a.s.
$$
In other words,
$$
Q_\omega\in\cap_{(p,q)\in\mI_{p_0}, T>\tau(\omega)}\cK^{p,q}_{\tau(\omega), T},
$$
which together with \eqref{DA0} yields that there is a $\mP$-null set $N$ such that for all $\omega\notin N$,
$$
Q_\omega\in\sC(\tau(\omega),\omega_{\tau(\omega)}).
$$

{\bf Verification of (C3)} Let $\mP\in\sC(s,x)$ and $\tau\geq s$ be a finite stopping time. For any $\cB_\tau$-measurable kernel 
$\mC\ni\omega\mapsto Q_\omega\in\cP(\mC)$ 
with 
$$
Q_\omega\in\sC(\tau(\omega),\omega_{\tau(\omega)}),\ \ \forall\omega\in\mC.
$$ 
By \cite[Theorem 6.1.2]{St-Va}, one knows that
\begin{align}\label{DA1}
\mP\otimes_\tau Q\in\cM^{\sigma,b}_{s,x}.
\end{align}
For fixed $p,q\in\mI^d_{p_0}$ and $T>s$, we want to show that there are $\kappa,\theta$ independent of $(s,x)$ 
such that for any $s\leq t_0<t_1\leq T$,
\begin{align}\label{DS0}
\mE^{\mP\otimes_\tau Q}\left(\int^{t_1}_{t_0} f(r,\omega_r)\dif r\Big|\cB_{t_0}\right)
\leq \kappa(t_1-t_0)^\theta\| f\|_{\mL^{q,p}_{t,x}},\ \mP\otimes_\tau Q-a.s.,
\end{align}
which means that
$$
\mP\otimes_\tau Q\in \cap_{(p,q)\in\mI_{p_0}, T>s}\cK^{p,q}_{s,T}.
$$
We make the following decomposition:
\begin{align*}
\mE^{\mP\otimes_\tau Q}\left(\int^{t_1}_{t_0} f(r,\omega_r)\dif r\Big|\cB_{t_0}\right)=I_1+I_2+I_3+I_4,
\end{align*}
where
\begin{align*}
I_1&:=\1_{\{\tau\leq t_0\}}\mE^{\mP\otimes_\tau Q}\left(\int^{t_1}_{t_0} f(r,\omega_r)\dif r\Big|\cB_{t_0}\right),\\
I_2&:=\1_{\{t_0<\tau\leq t_1\}}\mE^{\mP\otimes_\tau Q}\left(\int^{\tau}_{t_0} f(r,\omega_r)\dif r\Big|\cB_{t_0}\right),\\
I_3&:=\1_{\{t_0<\tau\leq t_1\}}\mE^{\mP\otimes_\tau Q}\left(\int^{t_1}_{\tau} f(r,\omega_r)\dif r\Big|\cB_{t_0}\right),\\
I_4&:=\1_{\{t_1<\tau\}}\mE^{\mP\otimes_\tau Q}\left(\int^{t_1}_{t_0} f(r,\omega_r)\dif r\Big|\cB_{t_0}\right).
\end{align*}
For $I_1$, noting that 
$$
I_1=\1_{\{\tau\leq t_0\}}\mE^{\mP\otimes_\tau Q}\left(\int^{t_1\vee\tau}_{t_0\vee\tau} f(r,\omega_r)\dif r\Big|\cB_{t_0\vee\tau}\right),
$$
by Lemma \ref{Le73} below, there is a $\mP\otimes_\tau Q$-null set $N\in\cB_\tau$ so that for all $\omega\notin N$,
\begin{align*}
\mE^{\mP\otimes_\tau Q}\left(\int^{t_1\vee\tau}_{t_0\vee\tau} f(r,\omega_r)\dif r\Big|\cB_{t_0\vee\tau}\right)
&=\mE^{Q_\omega}\left(\int^{t_1\vee\tau}_{t_0\vee\tau} f(r,\omega_r)\dif r\Big|\cB_{t_0\vee\tau}\right)\\
&\leq \kappa(t_1-t_0)^\theta\| f\|_{\mL^{q,p}_{t,x}},\ Q_\omega-a.s.
\end{align*}
Hence,
$$
I_1\leq \kappa(t_1-t_0)^\theta\| f\|_{\mL^{q,p}_{t,x}},\ \mP\otimes_\tau Q-a.s.
$$
For $I_2$, since $\mP\otimes_\tau Q|_{\cB_\tau}=\mP|_{\cB_\tau}$, we have
\begin{align*}
I_2&=\1_{\{t_0<\tau\leq t_1\}}\mE^{\mP}\left(\int^{\tau}_{t_0} f(r,\omega_r)\dif r\Big|\cB_{t_0}\right)
\leq \kappa(t_1-t_0)^\theta\| f\|_{\mL^{q,p}_{t,x}}.
\end{align*}
For $I_3$, since $(\mP\otimes_\tau Q)(\cdot|\cB_\tau)(\omega)=Q_\omega$, we have
\begin{align*}
I_3&=\1_{\{t_0<\tau\leq t_1\}}\mE^{\mP\otimes_\tau Q}\left(\int^{t_1}_{\tau} f(r,\omega_r)\dif r\Big|\cB_{t_0\wedge\tau}\right)\\
&=\1_{\{t_0<\tau\leq t_1\}}\mE^{\mP}\left(\mE^{Q_\cdot}\left(\int^{\tau}_{t_0\wedge\tau} f(r,\omega_r)\dif r\right)\Big|\cB_{t_0\wedge\tau}\right)\\
&\leq \kappa(t_1-t_0)^\theta\| f\|_{\mL^{q,p}_{t,x}}.
\end{align*}
Lastly, for $I_4$ we have
\begin{align*}
I_4&=\1_{\{t_1<\tau\}}\mE^{\mP}\left(\int^{t_1}_{t_0} f(r,\omega_r)\dif r\Big|\cB_{t_0}\right)
\leq \kappa(t_1-t_0)^\theta\| f\|_{\mL^{q,p}_{t,x}}.
\end{align*}
Combining the above calculations, we obtain \eqref{DS0}.
The proof is competed by Theorem \ref{Th72} below.
\end{proof}

\section{Examples}

For $R\geq 1$, let $\phi_R:[0,\infty)\to[0,\infty)$ be a smooth increasing function with
$$
\phi_R(r)=r,\ \ \ r\leq R;\ \ \phi_R(r)=R+1,\ \  r\geq 2R.
$$
For $\alpha\in\mR$ and $n\in\mN$, define
$$
f_R^{(\alpha)}(r):=(\phi_R(r))^\alpha,\ \ f_{R,n}^{(\alpha)}(r):=(\phi_R(r+\tfrac1n))^\alpha.
$$
Clearly,
$$
f^{(\alpha)}_R(r)=r^\alpha\ \mbox{ for }r<R\ \mbox{ and }
\ \ f_{R,n}^{(\alpha)}(r)=(r+\tfrac1n)^\alpha\ \mbox{ for }r+\tfrac1n<R.
$$
Below we provide two examples to illustrate the assumption {\bf ($\widetilde {\bf H}^\sigma$)}.
\bx\label{Ex51}\rm
Let $d\geq 3$ and $0<\alpha<(\frac d2-1)\wedge(\frac12+\frac{1}{d-1})$.
Let 
$$
\sigma(x)=f_R^{(-\frac\alpha2)}(|x|^2)\mI_{d\times d}.
$$ 
We verify {\bf ($\widetilde {\bf H}^\sigma$)} for $\sigma_n(x)=f_{R,n}^{(-\frac\alpha2)}(|x|^2)\mI_{d\times d}$.
Note that
$$
a_n(x)=(\sigma_n\sigma^*_n)(x)=f_{R,n}^{(-\alpha)}(|x|^2)\mI_{d\times d}.
$$
Thus,
$$
\lambda_n(x)=\mu_n(x)=f_{R,n}^{(-\alpha)}(|x|^2).
$$
In particular, we have
$$
\lambda^{-1}_n(x)\leq\phi^{\alpha}_R(|x|^2+1)\in \widetilde\mL^{\infty}(\mR^d),
$$
and for $p_1<\frac d{2\alpha}$,
$$
\mu_n(x)\leq \phi^{-\alpha}_R(|x|^2)\in \widetilde\mL^{p_1}(\mR^d).
$$
On the other hand, by the chain rule, we have
$$
\p_ia^{ij}_n(x)=2x_j(f^{(-\alpha)}_{R,n})'(|x|^2)
$$
and
$$
\p_i\p_ja^{jj}_n(x)=\Delta f^{(-\alpha)}_{R,n}(|\cdot|^2)(x)=2d (f^{(-\alpha)}_{R,n})'( |x|^2)+4|x|^2(f^{(-\alpha)}_{R,n})''(|x|^2).
$$
Note that
$$
(f^{(-\alpha)}_{R,n})'(r)=-\alpha \phi_R(r+\tfrac1n)^{-\alpha-1}\phi'_R(r+\tfrac1n)
$$
and
\begin{align*}
(f^{(-\alpha)}_{R,n})''(r)&=-\alpha \phi_R(r+\tfrac1n)^{-\alpha-1}(\phi'_R(r+\tfrac1n)+\phi''_R(r+\tfrac1n))\\
&\quad+\alpha(\alpha+1)\phi_R(r+\tfrac1n)^{-\alpha-2}(\phi'_R(r+\tfrac1n))^2.
\end{align*}
It is easy to see that for $p_1<\frac{d}{2\alpha+1}$,
$$
|\p_ia^{ij}_n(x)|\leq 2\alpha|x|^{-2\alpha-1}\1_{\{|x|^2\leq 2R\}}\in\widetilde\mL^{p_1}(\mR^d),
$$
and due to $\alpha<\frac d2-1$,
$$
\p_i\p_ja^{jj}_n(x)\leq C_{\alpha,R}.
$$
Hence, \eqref{NA0} holds for $p_0=\infty$, $p_1\in(\frac{d-1}{2},\frac{d}{2\alpha+1})$ and $p_2=q_2=\infty$.
Moreover, if $p_4<\frac{d}{\alpha}$, then
$$
|\sigma_n(x)|\leq \phi^{-\frac{\alpha}2}_R(|x|^2)\in \widetilde\mL^{p_4}(\mR^d).
$$
Thus, \eqref{NA1} holds for $p_4\in(d,\frac{d}{\alpha})$ and $q_4=\infty$. Therefore, {\bf ($\widetilde {\bf H}^\sigma$)} 
is satisfied for the above $\sigma_n(x)$.
In particular, by Theorem \ref{Th43}, there exists at least one solution for the following singular SDE:
$$
\dif X_t=\phi_R(|X_t|^2)^{-\alpha/2}\dif W_t+b(X_t)\dif t,\ \ X_0=x,
$$
where $\alpha\in(0,(\frac d2-1)\wedge(\frac12+\frac{1}{d-1}))$ and $b\in\widetilde L^{p}$ for some $p>\frac{d}{2}$ satisfies $(\div b)^-=0$.
\bp\label{Prop62}
Let $d\geq 3$, $\alpha\in(0,(\frac d2-1)\wedge(\frac12+\frac{1}{d-1}))$, $\beta\in(0,2\alpha)$ and $\lambda\geq 0$.
For each $x\in\mR^d$, the following SDEs admits a unique strong solution:
\begin{align}\label{SDE9}
\dif X_t=|X_t|^{-\alpha}\dif W_t+\lambda X_t|X_t|^{-\beta-1}\dif t,\ \  X_0=x.
\end{align}
\ep
\begin{proof}
Let $b(x):=\lambda x|x|^{-\beta-1}$, and for $R\in\mN$, let $\sigma_R(x)=\phi_R(|x|^2)^{-\alpha/2}\mI$. 
Since $\lambda\geq 0$ and $\beta<2$, it is easy to see that $b\in\widetilde L^{p}$ for any $p\in(\frac d2, \frac d\beta)$ and $(\div b)^-\equiv 0$.
Let $X^R_t$ solve the following SDE:
$$
X^R_t=x+\int^t_0\sigma_R(X_s^R)\dif W_s+\int^t_0b(X^R_s)\dif s.
$$
Let $\Phi:\mR_+\to\mR_+$ be a smooth function with $\Phi(r)=1$ for $|r|\leq 1$ 
and $\Phi(r)=r$ for $r>2$. By It\^o's formula, it is easy to see that
$$
\sup_{R\in\mN}\bE\left(\sup_{t\in[0,T]}\Phi(|X^R_t|^2)\right)\leq C.
$$
From this, by Chebyshev's inequality, we derive that
$$
\lim_{R\to\infty}\bP\left(\sup_{t\in[s,T]}|X^R_t|>R\right)=0,
$$
which together with Theorem \ref{Th55} implies that the assumptions of Theorem \ref{Th75} is satisfied.
So, there exists a solution to SDE \eqref{SDE9}.
To show the pathwise uniqueness, note that
$$
|\nabla b(x)|\leq C|x|^{-\beta-1},
$$
and for any $R>0$,
$$
\int_{B_R}|x|^{-(\beta+1) d}\det(\sigma\sigma^*)^{-1}(x)\dif x=\int_{B_R}|x|^{-(\beta+1) d+2\alpha d}\dif x<\infty.
$$
Thus by \cite[Theorem 1.1]{Wa-Zh} and the computations in Example 1 of \cite{Wa-Zh}, we obtain the uniqueness.
\end{proof}
\ex

\bx\rm
Let $d=2$ and $\alpha\in(0,\frac14)$. Consider the following diffusion matrix:
$$
\sigma(x)=\left(
\begin{aligned}
&f^{(\frac\alpha2)}_R(|x_2|^{2}),& 0\\
&0,&f^{(\frac\alpha2)}_R(|x_1|^{2})\\
\end{aligned}
\right).
$$
Let us define
$$
\sigma_n(x):=\left(
\begin{aligned}
&f^{(\frac\alpha2)}_{R,n}(|x_2|^2),& 0\\
&0,&f^{(\frac\alpha2)}_{R,n}(|x_1|^2)\\
\end{aligned}
\right),\ 
a_n(x):=\left(
\begin{aligned}
&f^{(\alpha)}_{R,n}(|x_2|^2),& 0\\
&0,&f^{(\alpha)}_{R,n}(|x_1|^2)\\
\end{aligned}
\right).
$$
Then
$$
\lambda_n(x)=f^{(\alpha)}_{R,n}(|x_2|^2)\wedge f^{(\alpha)}_{R,n}(|x_1|^2),\quad \mu_n(x)=
f^{(\alpha)}_{R,n}(|x_2|^2)\vee f^{(\alpha)}_{R,n}(|x_1|^2).
$$
Clearly, we have
$$
\lambda^{-1}_n(x)=f^{(-\alpha)}_{R,n}(|x_2|^2)\vee f^{(-\alpha)}_{R,n}(|x_1|^2)
$$
and
$$
\p_ia^{ij}_n(x)=\p_i\p_j a^{ij}_n(x)=0.
$$
Thus, \eqref{NA0} holds for any $p_0\in(2,\frac{1}{2\alpha})$, $p_1=\infty$ and $p_2=q_2=\infty$.
Moreover, it is easy to see that \eqref{NA1} holds for $q_4=\infty$ and any $p_4\in(\frac{2p_0}{p_0-1},\infty)$.
Therefore, {\bf ($\widetilde {\bf H}^\sigma$)} holds for the above $\sigma_n(x)$.
As in Proposition \ref{Prop62}, by Theorems \ref{Th55} and \ref{Th75}, 
for any starting point $X_0=x\in\mR^2$, there exists at least one solution for 
the following two dimensional degenerate SDE:
$$
\left\{
\begin{aligned}
\dif X^1_t=|X^2_t|^\alpha\dif W^1_t+b^1(X_t)\dif t,\\
\dif X^2_t=|X^1_t|^\alpha\dif W^2_t+b^2(X_t)\dif t,
\end{aligned}
\right.
$$
where $\alpha\in(0,\frac14)$ and $b=(b^1, b^2)\in\widetilde L^{p}(\mR^2)$ for some $p>\frac{2}{1-4\alpha}$, 
and for some $K\in\mN$,
$$
|b(x)|\leq C|x|,\ \ |x|>K.
$$
We would like to say some words about the range of $p$. Intuitively, when $X_t$ moves to the unit ball, smaller $\alpha$ means stronger noise and so the drift $b$ could be more singular.
While, the uniqueness for the above example is left open, even for $b=0$.
\ex

\section{Appendix}


We first recall the following lemma (cf. \cite[Theorem 6.1.2]{St-Va}).
\bl\label{Le71}
Let $\tau$ be a finite stopping time and $\mC\ni\omega\mapsto Q_\omega\in\cP(\mC)$ be a $\cB_\tau$-measurable probability kernel. 
Given a probability measure $\mP\in\cP(\mC)$, there exists a unique probability measure $\mP\otimes_\tau Q\in\cP(\mC)$ so that
$$
(\mP\otimes_\tau Q)|_{\cB_\tau}=\mP|_{\cB_\tau},\ \ (\mP\otimes_\tau Q)(\cdot|\cB_\tau)(\omega)=Q_\omega(\cdot).
$$
In particular,
$$
(\mP\otimes_\tau Q)(\Gamma)=\int_\mC Q_\omega(\Gamma)\mP(\dif\omega),\ \ \forall\Gamma\in\cB:=\vee_{t\geq 0}\cB_t.
$$
\el
For each $(s,x)\in\mR_+\times\mR^d$, let $\sC(s,x)$ be a non-empty convex subset of $\cP(\mC)$ with 
$$
\mP\{\omega: \omega_s=x\}=1.
$$
We suppose that $\{\sC(s,x): (s,x)\in\mR_+\times\mR^d\}$ satisfies
\begin{enumerate}[{\bf (C1)}]
\item  Let $(s_n,x_n)$ converge to $(s,x)$. For any sequence $\mP_n\in\sC(s_n,x_n)$, there is a subsequence $n_k$ and $\mP\in\sC(s,x)$ so that 
$\mP_{n_k}$ converges to $\mP$.
\item (Disintegration) Let $\mP\in\sC(s,x)$ and $\tau\geq s$ be a finite stopping time. For any r.c.p.d. $(\mP_\omega)_{\omega\in\mC}$ 
of $\mE^\mP(\cdot|\cB_\tau)$, there is a $\mP$-null set $N\in\cB_\tau$ so that
$$
\mP_\omega\in\sC(\tau(\omega),\omega_{\tau(\omega)}), \ \omega\notin N.
$$
\item (Reconstruction) Let $\mP\in\sC(s,x)$ and $\tau\geq s$ be a finite stopping time. For any $\cB_\tau$-measurable kernel 
$\mC\ni\omega\mapsto Q_\omega\in\cP(\mC)$ with 
$$
Q_\omega\in\sC(\tau(\omega),\omega_{\tau(\omega)}), \forall\omega\in\mC,
$$ 
it holds that
$$
\mP\otimes_{\tau}Q\in\sC(s,x).
$$
\end{enumerate}

We have the following strong Markov selection theorem, whose proofs are completely the same as in \cite[Theorem 12.2.3]{St-Va}
(see also \cite{Fl-Ro} and \cite[Theorem 2.7]{Go-Ro-Zh}). We omit the details.
\bt\label{Th72}
Under {\bf (C1)}, {\bf (C2)} and {\bf (C3)}, there is a measurable selection 
$$
\mR_+\times\mR^d\ni (s,x)\mapsto \mP_{s,x}\in\sC(s,x)
$$ 
so that for any 
$(s,x)\in\mR_+\times\mR^d$ and finite stopping time $\tau\geq s$,
$\omega\mapsto \mP_{\tau(\omega), x(\tau(\omega),\omega)}$ is a r.c.p.d. of $\mP_{s,x}$ with respect to $\cB_\tau$. 
More precisely, there is a $\mP_{s,x}$-null set $N\in\cB_\tau$
such that for all $\omega\notin N$,
$$
\mP_{s,x}(\cdot|\cB_\tau)(\omega)=\mP_{\tau(\omega), \omega_{\tau(\omega)}}(\cdot).
$$
\et

The following two simple lemmas are used in the proof of Theorem \ref{Th55} (see \cite{Go-Ro-Zh}).
\bl\label{Le73}
Let $\sG\subset\sC$ be two countably generated sub $\sigma$-algebras  of $\cB$. Given $\mP\in\cP(\mC)$, 
let $Q_\omega$ be a r.c.p.d. of $\mP$
with respect to $\sG$. Then there is a $\mP$-null set $N\in\sG$ depending on $\sC$ and $\xi$ such that for all $\omega\notin N$,
$$
\mE^\mP(\xi|\sC)=\mE^{Q_\omega}(\xi|\sC),\ \ Q_\omega-a.s.
$$
\el
\begin{proof}
Let $A\in\sG$ and $B\in\sC$. By definition, we have
\begin{align*}
\int_A\mE^{Q_\omega}(1_B\mE^\mP(\xi|\sC))\mP(\dif \omega)
&=\int_A\mE^\mP(1_B\mE^\mP(\xi|\sC)|\sG)(\omega)\mP(\dif \omega)\\
&=\mE^\mP(1_A1_B\xi)=\int_A\mE^{Q_\omega}(1_B\xi)\mP(\dif\omega)\\
&=\int_A\mE^{Q_\omega}(1_B\mE^{Q_\omega}(\xi|\sC))\mP(\dif\omega).
\end{align*}
Hence, for each $B\in\sC$, there is a $\mP$-null set $N_B\in\sG$ so that for all $\omega\notin N$,
$$
\mE^{Q_\omega}(1_B\mE^\mP(\xi|\sC))=\mE^{Q_\omega}(1_B\mE^{Q_\omega}(\xi|\sC)).
$$
Since $\sC$ is countably generated, one can find a common null set $N_{\xi,\sC}$ so that for all $\omega\notin N$ and $B\in\sC$,
$$
\mE^{Q_\omega}(1_B\mE^\mP(\xi|\sC))=\mE^{Q_\omega}(1_B\mE^{Q_\omega}(\xi|\sC)),
$$
which in turn yields the desired result.
\end{proof}
\bl\label{Le74}
Let $\tau$ be a finite stopping time and $Q_\omega$ be a r.c.p.d. of $\mP$ with respect to $\cB_\tau$.
Let $X_t$ be a bounded continuous process. Suppose that for any $t\geq 0$,
$$
\mE^\mP(X_t|\cB_t)\leq A,\ \ \mP-a.s.
$$
Then there is a $\mP$-null set $N\in\cB_\tau$ such that for all $\omega\notin N$ and $t\geq\tau(\omega)$,
$$
\mE^{Q_\omega}(X_t|\cB_t)\leq A,\ \ Q_\omega-a.s.
$$
\el
\begin{proof}
By Lemma \ref{Le44}, we have
$$
\mE^\mP(X_{t\vee\tau}|\cB_{t\vee\tau})\leq A,\ \ t\geq 0.
$$
By Lemma \ref{Le73}, there is a $\mP$-null set $N$ such that for all $\omega\notin N$ and all rational number $t>0$,
$$
\mE^{Q_\omega}(X_{t\vee\tau}|\cB_{t\vee\tau})\leq A,\ \ Q_\omega-a.s.
$$
For fixed $\omega\notin N$, since (cf. \cite[p34. (3.15)]{St-Va})
$$
Q_\omega\{\omega': \tau(\omega')=\tau(\omega)\}=1,
$$
we have for all rational number $t\geq\tau(\omega)$,
$$
\mE^{Q_\omega}(X_{t}|\cB_{t})=\mE^{Q_\omega}(X_{t\vee\tau}|\cB_{t\vee\tau})\leq A,\ \ Q_\omega-a.s.
$$
Now for general $t\geq\tau(\omega)$, let $t_n\downarrow t$ be rational numbers.
By the dominated convergence theorem, we have
$$
\mE^{Q_\omega}(X_t|\cB_t)=\lim_{t_n\downarrow t}\mE^{Q_\omega}(X_{t_n}|\cB_{t})
=\lim_{t_n\downarrow t}\mE^{Q_\omega}(X_{t_n}|\cB_{t_n}|\cB_t)\leq A.
$$
The proof is complete.
\end{proof}

The following result provides a way of constructing a global solution from local solutions.
\bt\label{Th75}
Suppose that for each $R\in\mN$ and $(s,x)\in\mR_+\times\mR^d$, 
there is at least one local martingale solution $\mP^R_{s,x}\in\cM^{\sigma_R,b_R}_{s,x}$ so that $(s,x)\mapsto\mP^R_{s,x}$ is Borel measurable,
where 
$$
\sigma_R(t,x):=\sigma(t,\chi_R(x)x),\ \ b_R(t,x):=b(t,\chi_R(x)x),
$$
and
$$
\chi_R(x)=1,\ |x|\leq 2^{R-1},\ \ \chi_R(x)=0,\ |x|>2^{R}.
$$
Fix $(s_0,x_0)\in\mR_+\times\mR^d$.
If for each $T>s_0$ and any choice of
$\mP^R_{s_0,x_0}$ from $\cM^{\sigma_R,b_R}_{s_0,x_0}$,
\begin{align}\label{AS9}
\lim_{R\to\infty}\mP^R_{s_0,x_0}\left(\sup_{t\in[s,T]}|\omega_t|>R\right)=0,
\end{align}
then there is at least one local martingale solution $\mP\in\cM^{\sigma,b}_{s_0,x_0}$.
In particular, there is a global weak solution $(\frak{F}, X,W)$ for SDE \eqref{SDE0}.
\et
\begin{proof}
Without loss of generality, we assume $(s_0,x_0)=(0,0)$.  
Let $\tau_0=0$. We define a sequence of stopping times recursively by
$$
\tau_n:=\inf\{t>\tau_{n-1}: |\omega_t|>2^{n-1}\}=\inf\{t>0: |\omega_t|>2^{n-1}\},\ \ n\in\mN.
$$
Let $\mP^n_{s,x}\in\cM^{\sigma_n,b_n}_{s,x}$ be as in the assumptions. 
Define for $n\in\mN$,
$$
Q^{n}_\omega:=\mP^{n+1}_{\tau_{n}(\omega),\omega_{\tau_{n}(\omega)}},\ \omega\in\mC.
$$
Since $(s,x)\mapsto\mP^{n+1}_{s,x}$ is measurable, 
$\omega\mapsto Q^n_\omega$ is a $\cB_{\tau_n}$-measurable probability kernel on $\mC\times\cB$, 
i.e., for each $\Gamma\in\cB$, $\omega\mapsto Q^n_\omega(\Gamma)$ is $\cB_{\tau_n}$-measurable, and for each $\omega\in\mC$, $Q_\omega\in\cP(\mC)$.
Let $\widetilde\mP_1\in\cM^{\sigma_1,b_1}_{0,0}$. Define for $n\geq 2$,
$$
\widetilde\mP_{n+1}:=\widetilde\mP_1\otimes_{\tau_1}\mQ^1\otimes_{\tau_2}\cdots\otimes_{\tau_{n}}\mQ^n.
$$
By the construction and Lemma \ref{Le71}, one sees that
$$
\widetilde\mP_{n+1}|_{\cB_{\tau_{n}}}=(\widetilde\mP_{n}\otimes_{\tau_{n}}\mQ^n)|_{\cB_{\tau_{n}}}=\widetilde\mP_{n}|_{\cB_{\tau_{n}}},
$$
and by \cite[Theorem 1.2.10]{St-Va},
$$
\widetilde\mP_n\in\cM^{\sigma_n,b_n}_{0,0}.
$$
Moreover, by \eqref{AS9}, for each $T>0$,
$$
\lim_{n\to\infty}\widetilde\mP_n(\tau_n<T)=0.
$$
Finally, by \cite[Theorem 1.3.5]{St-Va}, there is a unique $\mP\in\cP(\mC)$ so that for each $n\in\mN$,
$$
\mP|_{\cB_{\tau_{n}}}=\widetilde\mP_{n}|_{\cB_{\tau_n}}.
$$
The proof is complete.
\end{proof}

{\bf Acknowledgement:} The author is very grateful to Mathias Sch\"affner  for pointing out an error in the earlier version.


\begin{thebibliography}{999}

\bibitem{An-De-Sl}Andres S., Deuschel J.-D. and Slowik M.: Invariance principle for the random conductance model in a degenerate ergodic environment. {\it Ann. Probab.}, (2015), no. 4, 1866-1891.

\bibitem{An-Ch-De-Sl}Andres S., Chiarini A., Deuschel J.-D. and Slowik M.: Invariance principle for random walks with time-dependent ergodic  degenerate weights. {\it Ann. Probab.} 46 (2018), no. 1, 302-336

\bibitem{Be-Sc}Bella P. and Sch\"affner M.: Local Boundedness and Harnack Inequality for Solutions of Linear Nonuniformly Elliptic Equations. 
{\it Comm. Pure Appl. Math.}, https://doi.org/10.1002/cpa.21876.

\bibitem{Be-Sc1}Bella P. and Sch\"affner M.: Non-uniformly parabolic equations and applications to the random conductance model.
arXiv:2009.11535.


\bibitem{De-Gi}De Giorgi E.: Sulla differenziabilit\`a e l' analiticit\`a delle estremali degli integrali multipli regolari. {\it Mem. Accad. Sci. Torino.}, 
ser 3a, 3, 25-43(1957).

\bibitem{Fl-Ro}Flandoli F. and Romito N.: Markov selections for the 3D stochastic Navier–Stokes equations, 
{\it Probab. Theory Related Fields} 140 (2008) 407-458.

\bibitem{Go-Ro-Zh}Goldys B., R\"ockner M. and Zhang X.:
Martingale solutions and Markov selections for stochastic partial differential equations.
{\it Stochastic Processes and their Applications}, 119 (2009) 1725-1764.

\bibitem{Ha-Li} Han Q. and Lin F.: {\it Elliptic partial differential equations.} Courant Institute of Mathematical Sciences. NewYork, 1997.


\bibitem{Ig-Ku-Ry}Ignatova M.,  Kukavica I. and Ryzhik L.: The Harnack inequality for second-order parabolic equations with divergence-free drifts of low regularity.
{\it Commun. in Partial Differ. Equa.}  VOL. 41, NO. 2, 208-226(2016).

\bibitem{Ka-Sh}Karatzas I. and Shreve S.E.: {\it Brownian motion and stochastic calculus}. Graduate Texts in Math., Springer-Verlag, 1988.

\bibitem{Kr73}Krylov N.V.: The selection of a Markov process from a Markov system of processes, 
and the construction of quasidiffusion processes, {\it Izv. Akad. Nauk SSSR Ser. Mat.} 37 (1973) 691-708.

\bibitem{Kr80} Krylov N. V.: {\it Controlled diffusion processes}. Translated from the Russian by A. B. Aries. Applications of Mathematics,14. 
Springer-Verlag, New York-Berlin, 1980. xii+308 pp.

\bibitem{Kr20}  Krylov N. V.: On time inhomogeneous stochastic It\^o equations with drift in $L_{d+1}$. arXiv:2005.08831v1.

\bibitem{Mo}Moser J.: On Harnack's theorem for elliptic differential equations. {\it Comm. Pure Appl. Math.} 14, 577-591(1961).

\bibitem{Na}Nash J.: Continuity of solutions of parabolic and elliptic equations. {\it Amer. J. Math.} 80, 931-954(1958).

\bibitem{Na-Ur} Nazarov A. and Ural'tseva N.N.: The Harnack inequality and related properties for 
solutions of elliptic and parabolic equations with divergence-free lower-order coefficients. 
{\it St. Petersburg Mathematical Journal, }23(1), (2012), 93-115.

\bibitem{Sk}Skorokhod A.V.: {\it Studies in the theory of random processes}. New York: Dover, 1982.

\bibitem{St-Va} Stroock D. W., Varadhan S. R. S.: {\it Multidimensional diffusion processes},
{Grundlehren der Mathematischen Wissenschaften}, {233}, {Springer-Verlag, Berlin-New York}, {1979}

\bibitem{Tr} Trudinger N. S.: On the regularity of generalized solutions of linear, non-uniformly elliptic equations. {\it Arch. Rational Mech. Anal.} 42, 50-62 (1971).

\bibitem{Wa-Zh}Wang Z. and Zhang X.: Existence and uniqueness of degenerate SDEs with H\"older diffusion and measurable drift.
{\it J. Math. Anal. Appl.}, 484 (2020) 123679.

\bibitem{Zh-Zh}Zhang X. and Zhao G.:  Stochastic Lagrangian path for Leray solution of 3D Navier-Stokes equations. 
{\it Comm. Math. Phys.}, volume 381, pages491-525(2021).

\bibitem{Zh-Zh1}Zhang X. and Zhao G.: Singular Brownian Diffusion Processes. {\it Communications in Mathematics and Statistics,} pp.1-49, 2018. 
\end{thebibliography}
\end{document}